\documentclass[14pt]{extarticle}
\usepackage{amssymb,amsfonts,amsthm,amsmath}
\usepackage{mathrsfs,pifont}
\usepackage[all]{xy}

\usepackage[top=1.2in, bottom=1.2in, left=0.5in,
right=0.5in]{geometry}

\usepackage{tikz}
\usetikzlibrary{matrix,arrows}

\usepackage{hyperref}
\usepackage[backrefs,msc-links]{amsrefs}

\newcommand{\C}{\mathbb{C}}
\newcommand{\Ct}{\mathbb{C}^\times}
\newcommand{\Q}{\mathbb{Q}}
\newcommand{\Z}{\mathbb{Z}}
\newcommand{\R}{\mathbb{R}}

\newcommand{\bT}{\mathsf{T}}
\newcommand{\bA}{\mathsf{A}}
\newcommand{\bff}{\mathsf{f}}
\newcommand{\bAb}{\overline{\bA}}
\newcommand{\bE}{\mathsf{E}}
\newcommand{\bC}{\mathsf{C}}
\newcommand{\bY}{\mathsf{Y}}
\newcommand{\bV}{\mathsf{V}}
\newcommand{\bZ}{\mathsf{Z}}
\newcommand{\bF}{\mathsf{F}}
\newcommand{\bR}{\mathsf{R}}
\newcommand{\bQ}{\mathsf{Q}}
\newcommand{\Gr}{\mathsf{Gr}}

\newcommand{\bj}{\mathsf{j}}

\newcommand{\tbT}{\widetilde{\bT}}
\newcommand{\bP}{\mathbb{P}}
\newcommand{\cK}{\mathscr{K}}
\newcommand{\cT}{\mathscr{T}}
\newcommand{\cL}{\mathscr{L}}
\newcommand{\cP}{\mathscr{P}}
\newcommand{\cU}{\mathscr{U}}
\newcommand{\cF}{\mathscr{F}}
\newcommand{\cM}{\mathscr{M}}
\newcommand{\tM}{\widetilde{\mathscr{M}}}
\newcommand{\cN}{\mathscr{N}}
\newcommand{\cG}{\mathscr{G}}
\newcommand{\fF}{\mathsf{Fl}}
\newcommand{\cV}{\mathscr{V}}
\newcommand{\cW}{\mathscr{W}}
\newcommand{\cE}{\mathscr{E}}
\newcommand{\cB}{\mathscr{B}}

\newcommand{\cC}{\mathscr{C}}

\newcommand{\bla}{\boldsymbol{\lambda}}
\newcommand{\bph}{\boldsymbol{\varphi}}

\newcommand{\be}{\mathbf{e}}
\newcommand{\bv}{\mathbf{v}}

\newcommand{\fC}{\mathfrak{C}}
\newcommand{\fK}{\mathfrak{K}}

\newcommand{\lan}{\left\langle} 
\newcommand{\ran}{\right\rangle} 

\newcommand{\fP}{\mathfrak{P}}

\newcommand{\fg}{\mathfrak{g}}
\newcommand{\fgh}{\widehat{\mathfrak{g}}}
\newcommand{\fh}{\mathfrak{h}}
\newcommand{\fn}{\mathfrak{n}}
\newcommand{\fu}{\mathfrak{u}}
\newcommand{\fb}{\mathfrak{b}}

\newcommand{\rd}{/\!\!/\!\!/\!\!/}
\newcommand{\rdd}{/\!\!/}

\newcommand{\vth}{\vartheta} 

\newcommand{\cO}{\mathscr{O}}
\newcommand{\Hd}{{H}^{\raisebox{0.5mm}{$\scriptscriptstyle \bullet$}}}

\newcommand{\tO}{\widehat{\mathscr{O}}}
\newcommand{\vir}{\textup{vir}}
\newcommand{\flop}{\textup{flop}}

\newcommand{\fAttr}{\Attr^f}
\newcommand{\tw}{\textup{tw}}

\newcommand{\rM}{M}

\newcommand{\bTg}{\bT_\textup{gauge}}
\newcommand{\Wg}{W_\textup{gauge}}
\newcommand{\Gg}{G_\textup{gauge}}

\DeclareMathOperator{\Stab}{Stab}

\DeclareMathOperator{\Hom}{Hom}
\DeclareMathOperator{\Ker}{Ker}
\DeclareMathOperator{\Img}{Im}
\DeclareMathOperator{\Coker}{Coker}
\DeclareMathOperator{\Aut}{Aut}

\DeclareMathOperator{\Lie}{Lie}

\DeclareMathOperator{\const}{const}
\DeclareMathOperator{\Ell}{Ell}

\DeclareMathOperator{\chr}{char}
\DeclareMathOperator{\pt}{pt}
\DeclareMathOperator{\cochar}{cochar}

\DeclareMathOperator{\rk}{rk}
\DeclareMathOperator{\Pic}{Pic}

\DeclareMathOperator{\ind}{ind}
\DeclareMathOperator{\Spec}{Spec}

\DeclareMathOperator{\spec}{Spec}
\DeclareMathOperator{\Attr}{Attr}
\DeclareMathOperator{\supp}{supp}
\DeclareMathOperator{\ev}{ev}
\DeclareMathOperator{\Vx}{\mathbf{V}}
\DeclareMathOperator{\tVx}{\widetilde{\Vx}}
\DeclareMathOperator{\QM}{\mathsf{QM}}
\DeclareMathOperator{\End}{End}
\DeclareMathOperator{\Thom}{Thom}
\DeclareMathOperator{\diag}{diag}
\DeclareMathOperator{\Mon}{Mon}
\DeclareMathOperator{\Sol}{Sol}

\newcommand{\Ld}{{\Lambda}^{\!\raisebox{0.5mm}{$\scriptscriptstyle
      \bullet$}}\!}

\newcommand{\tQM}{\widetilde{\QM}}
\newcommand{\xt}{t}

\newtheorem{Theorem}{Theorem}

\newtheorem{Lemma}{Lemma}[section]
\newtheorem{Proposition}[Lemma]{Proposition}
\newtheorem{Corollary}[Lemma]{Corollary}

\theoremstyle{definition}

\newcommand{\Mbar}{\overline{M}}

\begin{document}

\title{Elliptic stable envelopes} 
\author{Mina Aganagic and 
  Andrei Okounkov} 
\date{}
\maketitle

\center{\emph{\small To Igor Krichever, with gratitude for inspiration and friendship.}}

\abstract{We construct stable envelopes in equivariant elliptic
  cohomology of Nakajima quiver varieties. In particular, this 
gives an elliptic generalization of the results of \cite{MO1}. 
We apply them to the computation of the 
monodromy of $q$-difference equations arising in the enumerative 
K-theory of rational curves in Nakajima varieties, including the quantum 
Knizhnik-Zamolodchikov equations.} 

\setcounter{tocdepth}{2}
\tableofcontents

\section{Introduction} 

\subsection{Different levels of stable envelopes} 

\subsubsection{}

Given an action of a torus $\bA$ on an algebraic variety $X$, one 
can define the attracting correspondence 
\begin{equation}
\Attr = \left\{ (x,y) , \lim_{a\to 0} a\cdot x = y \right\} \subset 
X \times X^\bA\label{attr_corr} \,,
\end{equation}
where $0$ is a point at infinity of $\bA$, or more precisely a
fixed point of a certain toric compactification of $\bA$. 
Cycles of the form \eqref{attr_corr} appear often 
in geometry and its applications to mathematical physics. Elementary examples are\ Schubert cells 
in the Grassmannian $\Gr(k,n)$, or their conormals in $X=T^*
\Gr(k,n)$. 

A disadvantage of the cycles \eqref{attr_corr} is that they may become unstable 
against small perturbation of the action: torus orbits, like 
gradient lines etc.,  can break under specialization, and the 
attracting set becomes smaller. One can improve the cycles 
\eqref{attr_corr} by taking the limit of attracting cycles for 
a small perturbation. This will be an
$\Aut(X)^{\bA}$-invariant cycle supported on the \emph{full} 
attracting set, which is the set of pairs 
$(x,y)$ that belong to a chain of closures of $\bA$-orbits. 

In practice, it is much more useful to have a characterization of 
improved cycles in terms that refer only to the original, unperturbed 
action, and not to a small perturbation which may not be explicit,
may break some symmetries, or may not be available altogether. 
The goal of stable envelopes is to provide just that. 

\subsubsection{}
Stable envelopes work best when $X$ is an equivariant symplectic 
resolution \cite{Kal2} and the action of $\bA$ preserves the symplectic
form. In this case, among the deformations of $X$ one may also 
consider the noncommutative ones, that is, \emph{quantizations}, 
for which \eqref{attr_corr} becomes the parabolic induction 
functor. 

Nakajima varieties \cites{Nak1,Nak2} form the largest and richest
family of equivariant symplectic resolutions known to date. In 
\cite{MO1}, the authors use stable envelopes to construct geometric
actions of certain quantum groups, called Yangians $\bY(\fg)$, on the cohomology
of Nakajima varieties. This extends and generalizes earlier work of 
Nakajima \cite{Nak3}, Varagnolo \cite{Vara}, and others, in which certain
smaller algebras $\bY(\fg_\textup{Kac-Moody})\subset \bY(\fg)$ 
were made to act by 
an explicit assignment on generators. 

One of the main applications in \cite{MO1} is a
description of the quantum cohomology of all Nakajima varieties in 
terms of this Yangian action. This extends, in particular, earlier
results of \cites{OP1,MauObl} on quantum cohomology of the Hilbert schemes of 
points of ADE surfaces. This theory finds important 
applications in enumerative geometry of sheaves on threefolds. 

\subsubsection{}

It is natural to ask whether stable envelopes may be lifted to
equivariant K-theory classes on $X \times X^\bA$ enjoying similar
properties. Here, one wants to work $\bT$-equivariantly, 
where $\bT\subset \bA$ is a maximal torus in $\Aut(X)^{\bA}$.  In particular, 
$\bT$ scales the symplectic form by a nontrivial character which 
we denote $\hbar$. Note that a limit argument does \emph{not} produce 
a well-defined $\bT$-equivariant K-class, because there is no 
$\bT$-equivariant deformation of the action. 

Reflecting this, in K-theory stable envelopes acquire an additional 
parameter --- a fractional line bundle $\cL \in \Pic(X) \otimes_\Z
\R$, called the \emph{slope}, see \cites{Opcmi}. 
The dependence on $\cL$ is 
piecewise constant, with walls forming a certain locally finite 
$\Pic(X)$-periodic arrangement of rational hyperplanes in $\Pic(X) \otimes_\Z
\R$. The intricacies of this arrangement reflect, among other things, the 
intricacies of quantizations of $X$, especially over a field of prime 
characteristic, as in the work of Bezrukavnikov and his
collaborators, see e.g. \cites{BezF,BK,BezLo,BM,Kal1}. 

{}From the representation-theoretic viewpoint,
$$
\fh = \Pic(X)\otimes_\Z\textup{field}
$$
is the Cartan subalgebra of $\fg$ and the walls correspond to 
roots of the loop algebra $\fgh=\fg[u^{\pm 1}]$. The K-theoretic 
lift of the construction of the Yangian gives an action of
$\cU_\hbar(\fgh)$, which is a Hopf algebra deformation of 
the universal enveloping algebra $\cU(\fgh)$. 

In particular, 
quantum difference equations in the K-theory of Nakajima 
varieties have been determined in terms of this action, see \cite{OS}. 
As before, this has direct application to K-theoretic enumeration of 
sheaves on threefolds or, more precisely, to K-theoretic Donaldson-Thomas
theory. 

\subsubsection{}\label{sDHL}
One natural direction for further generalizations is to lift stable
envelopes to Fourier-Mukai functor, that is, to a $\bT$-equivariant
complex of coherent sheaves on $X \times X^\bA$. This is pursued in 
\cite{DHLMO} and one is hoping, in particular, to categorify the 
$\cU_\hbar(\fgh)$-action along these lines.  The dependence of 
stable envelopes on the slope $\cL$ remains the same piecewise 
constant dependence. 

In this paper, we go in a different direction, and construct stable 
envelopes in equivariant \emph{elliptic} cohomology over $\C$. Perhaps the 
most striking new feature of the theory is that the piecewise constant 
dependence on $\cL$ is replaced by a \emph{meromorphic} 
dependence on 
$$
z \in \Pic(X) \otimes_\Z E
$$
where $E=\Ct/q^\Z$ is the elliptic curve of the cohomology theory. 
The piecewise constant dependence is recovered in the limit when 
one of the periods $-\ln(q)$ goes to $+\infty$, so that 
$$
- \Re \, \frac{\ln z}{\ln q} \to \cL \,.
$$
It is, of course, well known that in such limit elliptic functions
have a piecewise analytic limit, as exemplified by 
\begin{equation}
\lim \left(- \Re \, \frac{\ln z}{\ln q}\right) \in (k,k+1) \quad 
\Rightarrow \quad 
\lim \frac{\vth(az)}{\vth(z)}  = a^{k+\frac12} \,,\label{limit_theta}
\end{equation}
where $\vth$ is the classical 
odd, that is, the one with $\vth(x^{-1})=-\vth(x)$, theta function. It is given 
explicitly by \eqref{theta_function}, from which
\eqref{limit_theta} is immediate. In particular, the poles of 
elliptic stable envelopes in $z$ form a certain refinement of the 
roots of $\fgh$. 

Our main result is the construction of elliptic stable envelopes for 
Nakajima varieties given in Theorem \ref{t1}. As a consequence, 
we lift the $\cU_\hbar(\fgh)$-action to an elliptic 
quantum group.

\subsection{Pole subtraction and monodromy}

Elliptic stable envelopes solve a certain connection problem 
for difference equations in $K_\bT(X)$, which we call the 
\emph{pole subtraction} problem. 

\subsubsection{}

Enumerative $K$-theory of rational curves in Nakajima varieties is a 
source of interesting and important linear difference equations 
with regular singularities for a
$K_\bT(X)$-valued function of 
$$
z\in \bZ = \Pic(X) \otimes \Ct \,, 
$$
see \cite{Opcmi}. The shifts
in these difference equations are $z\mapsto q^\cL z$, where 
$\cL\in \Pic(X)$. There are commuting difference equations that 
shift equivariant variables by $q^{\cochar \bT}$, among which the
shifts by $q^{\cochar \bA}$ give equations with regular
singularities. For brevity, we call these equations {\em quantum difference 
equations}. These flat difference connections include many 
important difference equations of mathematical physics, including 
the quantum Knizhnik-Zamolodchikov equations \cite{FrenResh}. 

A central question about linear difference, or differential, equations 
is their monodromy. In particular, in the more traditional theory of quantum 
groups, one links the monodromy of qKZ equations associated to 
quantum affine Lie algebras to representations of 
elliptic quantum groups, see e.g. \cites{EFK,EtMour,FelICM,FTV,FrenResh,Konno,Mour,Stok}. This is a difference analog of the 
description of the monodromy of the classical Knizhnik-Zamolodchikov 
equations given by Kohno and Drinfeld. 

It is not too 
surprising that our elliptic stable envelopes enter the monodromy 
computations. In fact, one can identify precisely the 
difference equation problem that elliptic stable envelopes solve. 

\subsubsection{}\label{s_sep_reg} 

By a theorem of Deligne \cite{Del}, a flat differential connection has 
regular singularities if it has regular singularities along
any curve intersecting the singular locus generically and
transversally. In stark contrast to this,  it is 
very easy for a difference equation to 
have regular singularities in each 
group of variables, but not jointly. Indeed, it suffices to look at the
function 
$$
f(z,a) = \exp \dfrac{\ln z \ln a}{\ln q} 
$$
which solves 
$$
f(qz,a) = a f(z,a) \,, \quad  f(z,qa) = z f(z,a) \,. 
$$
This is regular for $z\to 0, a\ne 0$ and also for $z\ne 0, a\to 0$, but not 
regular at the point $(z,a)=(0,0)$. 

The quantum difference equations have precisely this feature: they
have regular singularities in K\"ahler variables $z$ and equivariant 
variables $a\in\bA$, but not jointly. 
Rather than 
a bug, this will turn out to be a very important feature of the
theory. 

\subsubsection{} 

Let $(z,a)=(0,0)$ be an irregular point as above. More precisely, 
both $z$ and $a$ vary in certain toric varieties and a choice of 
a fixed point $(0,0)$ in the product of those varieties has the 
following geometric meaning. 

The difference equation in $z$ has regular singularities on
a toric compactification $\overline{\bZ}\supset \bZ$ given by the fan of 
ample cones of all flops of $X$.  So, a choice of the point $z=0$ is 
the choice of $X$ among all possible flops. Similarly, a choice of the 
point $a=0\in \overline{\bA}$ is a choice of attracting manifolds 
as in \eqref{attr_corr}. To this data, we associate elliptic 
stable envelopes and we show in Theorem \ref{t_ps} that they 
can be interpreted as the following connection matrices for 
quantum difference equations.

In a neighborhood 
of the point $(z,a)=(0,0)$ we have two kinds of solutions to our difference 
equations: there are $z$-solutions, which are holomorphic for $z\ne 0$ and 
meromorphic in $a$, and there are $a$-solutions, for which it is 
the other way around. Enumerative geometry naturally provides
a basis of $z$-solutions which we call \emph{vertex functions} or
vertices for short. 

A basis of $z$-solutions can be transformed, by a certain triangular 
$q$-periodic matrix $\fP$, to a basis of $a$-solutions. We call this transition 
matrix the \emph{pole subtraction matrix} because,  in principle, it can be
computed by quite literally subtracting poles, see Section \ref{s_ps}. In Theorem \ref{t_ps} we show that, with suitable normalization,
this pole subtraction matrix is given by the elliptic stable envelopes. 

\subsubsection{}

The resulting $a$-solutions are uniquely determined by their 
appropriately interpreted initial conditions at $a=0$. It is 
easy to identify those with the vertex functions, 
that is, $z$-solutions,  for the fixed locus $X^\bA$, with a certain shift of 
K\"ahler variables, see Proposition \ref{p_a0}. We get the 
following diagram of meromorphic isomorphisms 
\begin{equation}
  \xymatrix{
&
*+<10pt>[F-:<7pt>] {\textup{$a$-solutions for $X$}}
\ar[dr]^{\textup{initial conditions at $a=0$}}\\
*+<10pt>[F-:<7pt>]{\textup{vertices for $X$}} \ar[ur]^{\fP}&& 
*+<10pt>[F-:<7pt>]{\textup{vertices${}^\lhd$
 for $X^\bA$}}
\ar[ll]_{\Stab_{X,a=0,\textup{$q$-diff}}}
}\,,
\label{diag_trian} 
\end{equation}
in which we defined the bottom arrow so that it commutes and the 
superscript in 
``vertices${}^\lhd$'' indicates the shift of K\"ahler variables. 
Solutions of quantum difference equations naturally define a
sheaf on the product of $\Ell_\bT(X)$ with $\Pic(X)
\otimes_\Z E$. Theorem \ref{t_ps} specifies the identification 
of the meromorphic bottom map in \eqref{diag_trian} with 
elliptic stable envelopes. 

It follows at once that the monodromy of the difference 
equations in $a$ is given by elliptic $R$-matrices. As to the 
monodromy in K\"ahler variables, it fits into the following 
commutative square 
\begin{equation}
  \xymatrix{
*+<10pt>[F-:<7pt>]{\textup{vertices for $X$}} 
\ar[dd]^{\textup{Monodromy}}&&&
*+<10pt>[F-:<7pt>]{\textup{vertices${}^\lhd$
 for $X^\bA$}} 
\ar[lll]_{\Stab_{X,a=0,\textup{$q$-diff}}}
\ar[dd]_{\textup{Monodromy${}^\lhd$}}
\\
\\
*+<10pt>[F-:<7pt>]{\textup{vertices for $X_\flop$}} 
&&&
*+<10pt>[F-:<7pt>]{\textup{vertices${}^\lhd$
 for $X^\bA_\flop$}}
\ar[lll]_{\Stab_{X_\flop,a=0,\textup{$q$-diff}}}
}\,. 
\label{diag_square} 
\end{equation}
In particular, in the case when $\bA$ acts on the 
framing spaces of Nakajima quiver varieties as in 
Section \ref{s_tensor_prod}, this becomes an equation 
for the coproduct of the monodromy. Such equations 
play the decisive role in all known ways to compute 
the monodromy. 

\subsubsection{}
A more categorical way to talk about the diagrams above 
is the following\footnote{This was taught to us by Pavel Etingof.}. 
For a fixed quiver, let $\cC$ be the category
of $\cU_\hbar(\fgh)$-modules given by the 
equivariant K-theories of Nakajima varieties $X$. Extending 
the scalars, we can make it linear over $q$-periodic functions of 
equivariant and K\"ahler variables. The vertex functions 
define a functor 
\begin{equation}
\Sol: K_\bT(X) \mapsto \bV(X)\,, \label{funSol}
\end{equation}
from $\cU_\hbar(\fgh)$-modules to modules over $q$-periodic 
functions of $a$ and $z$ given by solutions of the quantum 
difference equations.
A tensor structure on $\cC$ comes from 
K-theoretic stable envelopes for the action of the framing 
torus as in Section \ref{s_tensor_prod} and \cite{MO1,OS}. 
To make \eqref{funSol} a fiber functor, we need to 
give it a tensor structure, and it follows from the diagram 
\eqref{diag_trian} that 
\begin{equation}
\fP^{-1} \circ \mathfrak{J} : \left(\Sol(X_1) \otimes \Sol(X_2) \right)^{\triangleleft}
\to \Sol(X_1 \otimes X_2) \label{twist1}
\end{equation}
is the required structure, where $\mathfrak{J}$ is the inverse of the 
isomorphism with initial conditions at $a=0$ in \eqref{diag_trian}. 
This operator is essentially the fundamental solution to quantum 
Knizhnik-Zamolodchikov equation and, up to normalizations, is
the fusion operator discussed in
\cite{EtMour,EtSchiff,EtVar,Jimboetal} 
and many other papers. It is given by a universal element in 
a completion of $\cU_\hbar(\fgh)\otimes \cU_\hbar(\fgh)$, and 
so one can do the twist \eqref{twist1} in stages, by first 
twisting the category $\cC$ by $\mathfrak{J}$ to another 
tensor category $\cC'$ as in \cite{EtMour, EtSchiff,Jimboetal} and 
leaving $\fP^{-1}$ to be the tensor structure on the resulting functor 
$$
\Sol' : \cC' \to \textup{solutions of the quantum difference
  equations} \,. 
$$
By construction, the braiding in the category $\cC'$ is the 
monodromy of qKZ, which Corollary \ref{corMonqKZ} below identifies, 
up to an explicit gauge transformation, with the elliptic 
$R$-matrix; see also the above cited papers for prior results in 
this direction.  Thus $\cC'$ may be identified with the category of modules
over the elliptic quantum group provided by elliptic cohomology of 
Nakajima quiver varieties. 

The square \eqref{diag_square} then means that the monodromy of 
the quantum difference equations gives a tensor isomorphism between 
the functors $\Sol'$ for $X$ and its flop $X$.  Since $\Sol'$ has 
essentially no tensor automorphisms, this is a very strong 
constraint on the monodromy, which will be more fully explored
in a separate paper. 

\subsubsection{}
In this paper, we work with $q$-difference equations originating 
in the $K$-theoretic counts of rational curves in a Nakajima variety
$X$ and we describe their monodromy in the language of elliptic 
cohomology of $X$, where $q$ is the modulus of the elliptic 
curve $E=\Ct/q^\Z$. One can specialize $q\to 1$, and thus connect
the monodromy of the quantum differential equation for $X$ to 
K-theoretic stable envelopes for $X$. 

A detailed and powerful link between the monodromy of 
the quantum differential equation for a symplectic resolution $X$ 
and the properties of the quantization of $X$ in characteristic $p\gg 0$ 
has been proposed by Bezrukavnikov and his collaborators.\footnote{It appears, no published account of these conjectures is 
available at the time of writing, hopefully this will change soon 
\cite{BezOk}.} Our results here allow to make a very substantial progress 
towards these conjectures. For quivers of finite type, this is 
directly related to a conjecture of V.~Toledano-Laredo which 
identifies the monodromy of the trigonometric Casimir connection 
for a Lie algebra $\fg$ with the quantum Weyl group in
$\cU_\hbar(\fgh)$. See \cite{GTL} for recent progress towards
that conjecture.

\subsection{Further directions}

Other areas where we expect elliptic stable envelopes to be very 
useful include:
\begin{itemize}
\item[---] correspondences of boundary conditions of supersymmetric gauge theories in three dimensions,
\item[---] knot theory and categorification.
\end{itemize}
We will return to these elsewhere, see e.g.\ 
\cite{AO2}, here we only sketch some salient aspects.

\subsubsection{}
In this paper, we work with $X$ which is a Nakajima variety or a hypertoric variety. Then, $X$ has a physical interpretation as a moduli space of Higgs vacua, or Higgs branch, of a supersymmetric gauge theory in a certain class. The relevant class of theories are 3d gauge theories with ${\cal N}=4$ supersymmetry, studied for example in \cite{IS, HW, GW, BDG}. In this context, vertex functions of $X$ have a physical interpretation as well: they are the supersymmetric partition functions of the gauge theory on ${\mathbb C} \times S^1$, with a choice of a vacuum state at infinity. The supersymmetric partition function is an appropriate supertrace in the Hilbert space of the theory on ${\mathbb C}$. (The manifold is a twisted product, where in going around the $S^1$ one twists ${\mathbb C}$ by $q$.) 

Gauge theories in this class have an important duality called \emph{3d mirror symmetry}, which is closely related to the \emph{symplectic duality} (see \cite{NakCoul} for review). 
The duality relates pairs of 3d theories, exchanging their equivariant
and K\"ahler parameters, and Higgs and Coulomb branches. The Higgs
branch of the mirror theory, which we will denote by $X^\vee$, is
expected to coincide with the Coulomb branch of the original theory,
and vice versa. See \cite{GW, NakCoul, BFN} for recent 
progress. 

Exchanging the roles of K\"ahler and equivariant parameters in Theorem
\ref{t_diff_eq}, we get a set of difference equations satisfied by
vertex functions of $X^\vee$. Recall that, by construction, vertices
are holomorphic in K\"ahler parameters. The physical content of Theorem
\ref{t_ps} is the correspondence
\begin{equation}
  \xymatrix{
*+<10pt>[F-:<7pt>]{\textup{vertices for $X$}} && 
*+<10pt>[F-:<7pt>]{\textup{vertices
 for $X^\vee$}}
\ar@{<->}[ll]_{\Stab}
}\,,
\label{Stab_} 
\end{equation}
between vertex functions of a pair of
dual theories. This diagram highlights how stable 
envelopes for $X$ and $X^\vee$ are fundamentally the 
same objects. In fact, there exists \cite{AO1} an elliptic class on the 
product of $X$ and $X^\vee$ that specializes to 
elliptic stable envelopes for both $X$ and $X^\vee$. 
Note it is essential for this to treat K\"ahler variable on 
the same footing as the equivariant variables in the 
definition of elliptic stable envelopes.

\subsubsection{}

Instead of working with gauge theories on ${\mathbb C}\times S^1$, one
can replace ${\mathbb C}$ with a disk $D$ with suitable conditions
imposed on the $T^2$ boundary. From this perspective, the elliptic
stable envelope may be interpreted as a certain operator defined in terms of the gauge theory on $L\times T^2$, where $L$ is an interval, with boundary conditions imposed on each end. Nothing depends on the size of the interval $L$, and taking it to be zero, we get in effect a 2d theory on $T^2$. The graded index of the Hilbert space of the 2d theory, or more precisely, its elliptic genus, computes the matrix elements of stable envelopes.

This suggests one should be able to promote the elliptic stable envelopes
to a functor between categories of boundary conditions for $X^\bA$
and $X$. For $3$-dimensional
gauge theories, one is interested in the category (in fact, 
2-category) of boundary conditions. Objects of this category, i.e. different boundary conditions, can be obtained by coupling the bulk 3d gauge theory to a 
${\cal N}=(2,2)$-supersymmetric theory on the 
2-dimensional boundary, see e.g. \cite{Kap} for a general 
discussion of such categories. 

This will be further explored in a companion paper, which will also discuss a common generalization of the 
elliptic stable envelopes and categorification of the K-theoretic 
stable envelopes mentioned above in Section \ref{sDHL}.

\subsubsection{}
Applications to knot theory arise in the special case of Nakajima
quivers based on ADE-type Dynkin diagrams. In this case, the
difference equations of Theorem \ref{t_diff_eq} include the quantum Kniznik-Zamolodchikov equations of \cite{FrenResh}.

The vertex functions in this context can be related to supersymmetric
partition functions of a variant of a certain 6d "little string
theory" of ADE-type, together with codimension 4 defects, see \cite{AH} for a review. 
The 6d theory string reduces to 6d ${\cal N}=(0,2)$ conformal
field theory in the point particle limit. The gauge theory of the
previous subsection is the theory on defects of the 6d theory. The
K\"ahler parameters of the quiver are moduli of the 6d theory; the
equivariant parameters describe positions of defects on a Riemann
surface ${\cal C} = {\mathbb C}^*$ on which 
6d theory is 
supported\footnote{
The full 6d geometry setting is very similar to that in
  \cite{AH}.}.  The elliptic $R$ matrices of section \ref{s_R_matr} describe braiding of defects on ${\cal C}$.

In the limit in which the qKZ equation reduces to the KZ equation, the little string theory reduces to the more familiar  6d CFT. The relation of this theory to knot theory and invariants of quantum groups has long been predicted by physicists \cite{OV, WF, GW2}. In a sense, this paper provides a significant step forward in realizing the physics prediction by providing a derivation of $R$-matrices from quantum field theory and string theory. It should be noted that to establish this result, one has to work with quantum K-theory, quantum affine algebras and the little string theory.

\subsection{Acknowledgments}
\label{sec:acknowledgments}

Our work on this paper was greatly facilitated by discussions with 
Andrey Smirnov, Tudor Dimofte, Davide Gaiotto,  Nikita Nekrasov, and Davesh Maulik. We are 
especially grateful to Nora Ganter for giving us the courage to work
in equivariant elliptic cohomology and 
to Pavel Etingof for his involvement and guidance on 
many occasions. 

AO thanks the Simons foundation for being financially 
supported as a Simons investigator and NSF for supporting 
enumerative geometry at Columbia as a part of FRG 1159416. 
 MA is supported, in part, by the NSF grant \#1521446.

As if by magic, many threads of our narrative, including elliptic
cohomology, integrable systems, and difference equations, all
converge in one person: Igor Krichever. Of course, the scientific explanation for
this apparent
magic is Igor's ability to see what is really deep and 
important long before others catch up. The writing and the
publication of this paper has been a long process indeed, perhaps the fate had
in mind that the paper should appear in time for Igor's 70th
birthday. We are very happy to use this chance to dedicate the paper
to Igor as a birthday present and thank him for the influence he had
on this paper, on us, and on the mathematical physics as a whole.

\section{Equivariant elliptic cohomology} 

\subsection{The curve $E$.} 

\subsubsection{} 

We set $E=\Ct/q^\Z$. This is a family of complex elliptic curves over
the punctured disc $0<|q|<1$. While in the general development of
elliptic cohomology it is very important to work with more general families of
elliptic curves, the above choice suffices for our purposes. Note, in
particular, that $E$ has no nontrivial endomorphisms and, more
generally, 
$$
\Hom(E^n,E^m) \cong \Hom(\Z^n,\Z^m) \,. 
$$

\subsubsection{} 

The theta function 
\begin{equation}
  \label{theta_function}
  \vth(x) = (x^{1/2} -x^{-1/2}) \, \prod_{n>0} (1-q^n x) (1-q^n/x) 
\end{equation}
satisfies
\begin{equation}
  \label{theta_transform}
  \vth(q^kx) = (-1)^k q^{-k^2/2} x^{-k} \vth(x) \,, \quad k\in \Z \,, 
\end{equation}
and thus defines a section of a degree $1$ line bundle on $E$ with the
unique zero at $x=1$. Translates of this line bundle form 
$\Pic_1(E)\cong E$, where $\Pic_1(E)$ denotes line bundles of
degree $1$. Similarly, meromorphic 
sections $s(x)$ of line bundles of degree
$d$ on $E$ satisfy
$$
s(qx) = c x^{-d} s(x)\,, \quad c \in \Ct \,, 
$$
see e.g.\  \cite{BikLang} for a systematic discusion of 
how line bundles on abelian varieties are described by their factors of
automorphy.

\subsubsection{} 

A long-standing tradition in the theory of elliptic functions is to
use additive notation in abstract formulas involving the group operation on abelian
varieties. At the same time, one uses $E=\Ct/q^\Z$ with its
multiplicative group law in concrete computations. We follow this
tradition, which is particularly convenient in the geometric
context. 

Indeed, in the world of complex-oriented cohomology theories,
the group law on $E$ comes from the \emph{tensor product} on line
bundles,  and so from the group operation in $GL(1)$. We believe the reader will appreciate the convenience of
translating vector bundles into elliptic functions with this choice of
notation.

\subsection{Basics}

\subsubsection{}

In this paper, we work with torus-equivariant elliptic cohomology 
over $\C$; this is sufficiently general to cover all applications that 
we have in mind. 

With one important exception in Section \ref{sfF}, our 
setting will be algebraic, that is, 
we will consider an algebraic torus
$\bT\cong (\C^\times)^n$ and regular equivariant maps $f: X\to Y$ between 
complex quasiprojective $\bT$-varieties.  

We also assume the $\bT$-action on $X$ is linearized which, 
by 
definition, means that the quasi-projective embedding of $X$ may 
be taken of the form 
\begin{equation}
X \hookrightarrow \bP(\textup{$\bT$-module})  \,.
\label{linearized}
\end{equation}
%

\subsubsection{}

Equivariant elliptic cohomology, developed in 
\cites{Groj,GKV,Rosu,Lurie,Gepner,Ganter} and other papers, defines a functor
$$
\Ell_\bT(X) : 
\big\{ \textup{$\bT$-spaces X} \big\} \to \{\textup{schemes}\}
$$
covariant in both $\bT$ and $X$, such that 
$$
\Ell_{(\Ct)^n}(\pt) \cong E^n
$$
for an elliptic curve $E$.  Strictly speaking, one should consider 
supercommutative schemes for varieties with odd cohomology, but 
Nakajima varieties and symplectic resolutions in general have only 
even cohomology \cite{Kal2}.

For a torus $\bT$, its characters and
  cocharacters are dual lattices defined by 
$$
\chr(\bT) = \Hom(\bT,\C^\times)\,, \quad \cochar(\bT) =
\Hom(\C^\times,\bT)\,.
$$
For the dual torus $\bT^\vee$, these are exchanged. 
Covariance in $\bT$ implies
$$
\Ell_\bT(\pt) = \bT / q^{\cochar(\bT)}=: \cE_\bT\,,
$$
canonically.

\subsubsection{}\label{s_small_analytic} 
To the projection $X \to \pt$ the functor $\Ell_\bT$ associates 
the map 
$$
\pi: \Ell_\bT(X) \to \cE_\bT
$$
that looks as follows in a small analytic neighborhood $U$ of 
a point $t\in \bT$: 
\begin{equation}
  \label{Ellt}
\xymatrix{ 
 \Spec H^*_{\bT}(X^{\bT_\xt},\C)  \ar[d] &  
\pi^{-1}(U) \ar@{->}[r] \ar[l] \ar[d]  & \Ell_\bT(X) \ar[d] \\
\Lie \bT& U \ar@{->}[r] 
\ar[l]_{\ln(\,\cdot\,)-\ln(t)} &\cE_\bT \,.
}
\end{equation}
Here both squares are pullbacks and the subgroup  
$$
\bT_\xt = \bigcap_{\chi(\xt)\in q^\Z} \Ker \chi \subset \bT \,,
$$
is the intersection of kernels of all characters 
$$
\chi \in \chr(\bT) \cong \Hom(\cE_\bT, E)  
$$
that are trivial on the image of $t$ in $\cE_\bT$. 

\subsection{Equivariant formality} 

\subsubsection{}

Let a torus $\bT$ act on a Nakajima variety $X$ so that is
scales its 
canonical symplectic form nontrivially. The fixed locus $X^\bT$ is then a smooth projective 
variety, which is a union of finitely many components 
$\{F_i\}$. A generic one-parameter subgroup 
$$
\Ct \owns s \mapsto \bT
$$
may be chosen so that the limit $\lim_{s\to 0} s \cdot x 
\in X^{\bT}$ exists for all $x\in X$.  For any 
subgroup $\bT'\subset \bT$, the terms in the 
corresponding 
Bia\l{}ynicki-Birula decomposition 
\begin{equation}
X^{\bT'} = \bigcup_{F_i} \, \big\{x \in X^{\bT'} \, \big| \, 
\lim_{s\to 0} s \cdot x 
\in F_i\big \} \label{BBx}
\end{equation}
are bundles of affine spaces over $F_i$ by \eqref{linearized}
and the classical 
result of \cite{BB}. The terms in the decomposition 
\eqref{BBx} are naturally partially ordered by containment in the closure. 

\subsubsection{}\label{s_flatness} 

The flatness of the left vertical arrow in \eqref{Ellt}, equivalently,
the freenees of $H^*_{\bT}(X^{\bT_\xt})$ over $H^*_{\bT}(\pt)$
is a property known as \emph{equivariant formality}, see  \cite{GKM}
for a comprehensive discussion. Formality implies that 
$$
H^*_{\bT}(X^{\bT_\xt}) \cong H^*(X^{\bT_\xt}) \otimes
H^*_{\bT}(\pt)
$$
and therefore 
$$
\pi^{-1}(t) \cong \Spec H^*(X^{\bT_\xt}) \,. 
$$

\begin{Lemma}
For any Nakajima variety $X$ and any subgroup $\bT'\subset \bT$, the fixed locus 
$X^{\bT'}$ is equivariantly formal. 
\end{Lemma}

\begin{proof}
By Corollary 1.3.2 
By part (8) of Theorem 14.1 in \cite{GKM}, it suffices to see that
the homology $H_*(X^{\bT'})$ is generated by 
$\bT$-invariant cycles. This follows from the 
Bia\l{}ynicki-Birula decomposition \eqref{BBx} by 
induction on the partial order and the long exact sequence of a pair. 
\end{proof}

\begin{Corollary}\label{cor1} 
For any subgroup $\bT'\subset \bT$, the cohomology $H^*(X^{\bT'})$ is even and 
$H^*(X^{\bT'},\Q)  \cong K^{\textup{top}}(X^{\bT'}) \otimes \Q
$. 
\end{Corollary}

\begin{proof}
By equivariant localization and formality, $H^*(X^{\bT})$ is a localization of 
$H^*_{\bT}(X) \cong H^*(X) \otimes
H^*_{\bT}(\pt)
$, which is even. The Bia\l{}ynicki-Birula decomposition then 
implies $H^*(X^{\bT'})$ is even for any subgroup $\bT'$. The 
comparison with the topological K-theory $K^{\textup{top}}(X^{\bT'})$ follows from 
the corresponding degeneration of the Atiyah-Hirzebruch spectral
sequence. 
\end{proof}

\subsection{Tautological generation} 

\subsubsection{}

In this revised version of the paper, we can take advantage of the
following powerful result of K.~McGerty and T.~Nevins which was not yet available at the time of the
writing. 

\begin{Theorem}[\cite{McGN}]\label{t_McGN} If $X$ is a Nakajima variety then 
$K^{\textup{alg}}_\bT(X) = K^{\textup{top}}_\bT(X)$ is generated by tautological 
bundles and $H^*_{\bT}(X,\Z)$ is generated by the Chern classes of the
tautological bundles. 
\end{Theorem}

As a corollary, $\Pic(X)$ and $\Pic_\bT(X)$ are lattices generated by tautological line
bundles.  They fit into an exact sequence 
\begin{equation}
0 \to \chr(\bT) \to \Pic_\bT(X) \to \Pic(X) \to 0 \label{PicXexact}
\end{equation}
of free abelian groups. 

\subsubsection{}

For completeness, we recall the proof of the following 
well-known fact, see e.g.\ \cite{CG, Merk}. 

\begin{Lemma}\label{lem_fib} 
In both algebraic and topological equivariant K-theory, we
have the following pullback diagram 
\begin{equation}
  \label{Kft}
\xymatrix{ 
 \Spec K(X^{t}) \ar[d] \ar[r] &  \Spec K_\bT(X)\ar[d] \\
\{t\} \ar@{->}[r] 
&\bT \,. 
}
\end{equation}
\end{Lemma}

\begin{proof}
Let $\chi$ be a character of $\bT$. Let $\C_\chi$ be $1$-dimensional $\bT$-module with 
character $\chi$ and 
$$
\Ct_\chi = \C_\chi \setminus \{0\} = \Img \chi \cong 
\bT/\Ker \chi \,. 
$$
We conclude 
\begin{align}
  \label{specializk}
 K_{\Ker \chi} (X)  &= K_\bT( X  \times \Ct_\chi ) \notag \\
 & =  \Coker \left(K_\bT(X) \xrightarrow{1-\chi} K_\bT(X) 
\right) \,, 
\end{align}
where the first line is 
Corollary 5 in \cite{Merk}, and
 the second line follows from the localization long
exact sequence applied to the inclusion 
$$
X \to X \times \C_\chi \,, 
$$
compare with Corollary 27 in \cite{Merk}. 

Since every subgroup $\bT' \subset \bT$ is a complete 
intersection kernels of characters, we have the 
following pullback diagram 
\begin{equation}
  \label{Kft1}
\xymatrix{ 
 \Spec K_{\bT'}(X) \ar[d] \ar[r] &  \Spec K_\bT(X)\ar[d] \\
\bT' \ar@{->}[r] 
&\bT \,. 
}
\end{equation}
On the other hand, by equivariant localization 
\begin{equation}
  \label{Kft2}
\xymatrix{ 
 U \times \Spec K(X^{\bT'}) \ar[d] \ar[r] &  \Spec K_{\bT'}(X)\ar[d] \\
U \ar@{->}[r] 
&\bT' \,, 
}
\end{equation}
where $U = \bT' \setminus \bigcup_{\bT'' \subsetneq \bT'} \bT''$ is
the open set of elements that generate a Zariski dense subgroup of
$\bT'$.  Applying \eqref{Kft1} and \eqref{Kft2} to the Zariski 
closed subgroup $\bT'$ generated by $t$, we obtain 
\eqref{Kft}. 
\end{proof}

As a corollary, $K^{\textup{alg}}_\bT(X^t) =
K^{\textup{top}}_\bT(X^t)$
for all $t\in \bT$, and we do not distinguish between 
these groups in what follows.

\subsection{Characteristic classes}

\subsubsection{} 

An equivariant rank $r$ complex vector bundle $V$ over $X$ defines a map
\begin{equation}
c: \Ell_\bT(X) \to \Ell_{GL(r)}(\pt) = S^r E \,, \label{defc}
\end{equation}
see Section (1.8) in \cite{GKV} and Section 5 in \cite{Ganter}. 
The coordinates in the target of \eqref{defc} are 
symmetric functions on $E^r$ --- symmetric functions in elliptic 
Chern roots.

\subsubsection{} 

Let $X$ be a Nakajima quiver variety, which is the case of main 
interest for us in this paper. By construction, $X$ is a quotient 
by $G= \prod GL(\bv_i)$. This gives a collection of 
tautological vector bundles $\{V_i\}$ of rank 
$\rk V_i = \bv_i$ and the map
\begin{equation}
  \label{HtoE}
  \Ell_\bT(X) \to \cE_\bT \times \prod S^{\bv_i} E \,. 
\end{equation}
By Corollary \ref{cor1} and Lemma \ref{lem_fib}, locally on 
$\cE_\bT$, this map may be modeled by the map 
\begin{equation}
  \label{KtoE}
  K_\bT(X) \to \bT \times \prod S^{\bv_i} \Ct 
\end{equation}
given by the K-theoretic Chern roots.

Theorem \ref{t_McGN} implies \eqref{KtoE} is an 
embedding, therefore \eqref{HtoE} is also an embedding.

\subsection{Pushforwards}

\subsubsection{}

Pullback in elliptic cohomology are the functorial maps 
$$
\Ell(f): \Ell_\bT(X) \to \Ell_\bT(Y)
$$
associated to a map $f: X \to Y$ of $\bT$-spaces. 
Pushforwards are defined for complex oriented maps 
and are sheaf homomorphisms
\begin{equation}
f_*:  \quad \Ell(f)_* \, \Theta(-N_f) \to \cO_{\Ell_\bT(Y)}\label{f_*}
\end{equation}
where $N_f \in K_\bT(X)$ is the normal bundle to $f$ and 
\begin{equation}
\Theta: K_\bT(X) \to \Pic \left(\Ell_\bT(X)\right) 
\label{Thom_map}
\end{equation}
is the \emph{Thom class} map defined as follows. 

\subsubsection{}

Let $V$ be a complex vector bundle over $X$.
Its Thom class is, by definition 
$$
\Theta(V) = c^* \cO(D) 
$$
where $c$ is the map \eqref{defc} and 
the divisor 
$$
D = \{0\} + S^{r-1} E \subset  S^{r} E
$$ is 
formed by those $r$-tuples that contain $0$. Since clearly 
$$
\Theta_{V_1 \oplus V_2}  = \Theta_{V_1} \otimes \Theta_{V_2}
$$
this extends to a group homomorphism \eqref{Thom_map}. 

\subsubsection{}

In English, the need to introduce the twist by the Thom class 
may be explained as follows. As we progress from 
equivariant cohomology to equivariant elliptic cohomology, 
the Euler class of the normal bundle becomes replaced by 
$\prod \vth (x_i)$, where 
$(x_1,\dots,x_r) \in S^r E$ are the Chern 
roots of the normal bundle. This product of $\vth$-functions is a section of 
a nontrivial line bundle over $S^r E$, which this very section
identifies with $\cO(D)$.

Also note that by construction of Chern classes, any bundle
$V$ of rank $r$,
whether it splits into line bundles or not, defines a map from 
$\Ell_\bT(X)$ to $S^r E$. The Thom class $\Theta(V)$ of $V$ is the pull-back of
$\prod \vth (x_i)$ under this map. In general, computations with Chern
roots of a vector bundle are computations with the coordinates $x_i$
on $S^r E$, which may be pulled back to $\Ell_\bT(X)$ via the Chern class
map.

\subsection{Universal line bundle $\cU$}

\subsubsection{}\label{s_tilde_c} 

For line bundles, the Chern class \eqref{defc} 
gives a group homomorphism 
$$
\Pic_\bT(X) \xrightarrow{\,\,c\,\,} \textup{Maps}\big(\Ell_\bT(X) \to E\big) 
$$
where the group operation in the target is 
the pointwise addition in $E$. This can be viewed as a map
$$
\tilde{c}: \Ell_\bT(X) \to  \cE_{\Pic_\bT(X)}^\vee \,,
$$
where
\begin{align}
\cE_{\Pic_\bT(X)} &= \Pic_\bT(X) \otimes_\Z E\,, \notag\\
  \cE_{\Pic_\bT(X)}^\vee &= \Hom(\Pic_\bT(X),E)  \label{Epic} 
\end{align}
is a pair of dual abelian varieties. The universal line bundle 
is a family of line bundles on $\Ell_\bT(X)$ pulled back via 
this map. 

\subsubsection{} 

On the product of two dual abelian varieties there is a universal
 line bundle $\cU_\textup{Poincar\'e}$. In complex analytic 
terms, the sections of $\cU_\textup{Poincar\'e}$ on  $E^\vee \times E$
are analytic functions on the universal cover with 
the same factors of automorphy as
\begin{equation}
\frac{\vth(sz)}{\vth(s) \, \vth(z)} \,, \label{cUthth}
\end{equation}
where we use the isomorphism $E^\vee \cong E$ given by the 
divisor  $\{1\}=(\vth) \subset E$. Note from \eqref{theta_transform}
that 
the function $\psi(z,s)$ as in \eqref{cUthth} 
satisfies
$$
\psi(qs,z)=z^{-1} \psi(s,z)\,, \quad
\psi(s,qz)=s^{-1} \psi(s,z)\,.
$$
Below, it will be convenient to consider pullbacks of $\cU_\textup{Poincar\'e}$ by 
the automorphism of the base. For example 
\begin{equation}
(z\mapsto z^{-1})^* \, \cU_\textup{Poincar\'e} \cong
\cU_\textup{Poincar\'e}^\vee\label{dualU} \,. 
\end{equation}

\subsubsection{} 
We define 
$$
\cU = (\tilde{c} \times 1 )^* \,  \cU_\textup{Poincar\'e} \,.
$$
This is a line bundle on 
\begin{equation}
\bE_\bT(X) = \Ell_\bT(X) \times
\cE_{\Pic_\bT(X)}\label{bET}
\end{equation}
which is a scheme over 
\begin{equation}
\cB_{\bT,X} = \cE_\bT \times
\cE_{\Pic_\bT(X)} \,.\label{bBT}
\end{equation}
We call the variables in the two factors of $\cB_{\bT,X}$ the
equivariant  and the K\"ahler parameters, respectively. 

\subsubsection{} 
The universal line bundle may be described very 
concretely using the map \eqref{HtoE}. The 
pullback under the addition map $S^k E \to E$ induces 
the isomorphisms 
\begin{equation}
\Pic_0 (S^k E) \cong \Pic_0(E) \cong E   \,.\label{PicSk}
\end{equation}
Meromorphic sections of a line bundle corresponding to 
$z\in E$ are symmetric meromorphic functions 
$\psi(s_1,\dots,s_k)$, $s_i
\in \Ct$, such that 
$$
\psi(q s_1, s_2, \dots,s_k) = z^{-1} \psi(s)  \,. 
$$
We get a factor of \eqref{PicSk} for each factor in 
the right-hand side of \eqref{HtoE}, including 
$\rk \bT$ many factors in $\cE_\bT$. This gives a map 
$$
\cE_{\Pic_\bT(X)}  \to \Pic_0(\cE_\bT \times \prod S^{\bv_i} E) \,, 
$$
and the bundle $\cU$ is pulled back from the corresponding 
bundle on the ambient space $\cE_{\Pic_\bT(X)}  \times \cE_\bT \times
\prod S^{\bv_i} E$. 

\subsubsection{}
Even more concretely, a Nakajima quiver variety 
is constructed as a GIT quotient a by
$$
\Gg=\prod GL(\bv_i)  \,,
$$
which is the complexified gauge group in the physical 
context. We have 
\begin{align}
\cE_\bT \times \prod S^{\bv_i} E &= 
\cE_{\bT \times \bTg} \big/ \Wg \,, \\
\cE_{\Pic_\bT(X)} & = \left(\cE^\vee_{\bT \times \bTg}\right)^{\Wg} \,, 
\end{align}
where 
$$
\bTg, \Wg \subset \Gg
$$
are the maximal torus and the Weyl group, respectively. 
Setting 
$$
\tbT=\bT \times \bTg\,,
$$
we have the 
natural map 
\begin{equation}
\left(\tbT^\vee\right)^{\Wg} \times \left(
\tbT\big/ \Wg \right)  \to \cE_{\Pic_\bT(X)}  \times \cE_\bT \times
\prod S^{\bv_i} E\label{coverbtT} \,.
\end{equation}
Meromorphic sections of the universal bundle pull back under
\eqref{coverbtT} to 
meromorphic functions of 
$z\in \left(\tbT^\vee\right)^{\Wg}$ and $s\in\tbT$ that 
are $\Wg$-invariant in $s$ and satisfy 
\begin{align}
\psi(z,q^\sigma s) &= \sigma(z)^{-1} \psi  \notag \,, \\
\psi(q^{\chi} z, s) &= \chi(s)^{-1} \psi \,, \label{trans} 
\end{align}
for any 
\begin{align*}
\sigma &\in \cochar \tbT = \chr \tbT^\vee\,, \\
\chi &\in \cochar \left(\tbT^\vee\right)^{\Wg} = 
\chr \left(\tbT \right)^{\Wg} \,. 
\end{align*}

\subsection{Shifts of K\"ahler variables}

\subsubsection{}

Observe that translations along the $\cE_{\Pic_\bT(X)}$ factor in 
\eqref{bET} preserve everything except the universal 
bundle $\cU$. Note that we may translate by an amount that 
depends on where we are in $\Ell_\bT(X)$. An example 
of such transformation is 
\begin{equation}
z \xrightarrow{\quad \tau(\lambda \mu) \quad}  z+ 
\lambda(\mu(t))\,, \label{def_tau}
\end{equation}
which may be defined for every pair
\begin{align*}
\mu &\in \chr(\bT) = \Hom(\cE_\bT,E)\,, \\
\lambda &\in \Pic_\bT(X) = \Hom(E, \cE_{\Pic_\bT(X)})  \,. 
\end{align*}
In \eqref{def_tau},  $t$ is a coordinate on $\cE_\bT$ and 
all coordinates in $\Ell_\bT(X)$ are unaffected by 
the transformation \eqref{def_tau}. 

\subsubsection{}
Let $\psi$ be a function on $\tbT^\vee \times \tbT$ as in 
\eqref{trans} and consider 
$$
\psi'(z,s) = \frac{\psi(\lambda(\mu(s)) z, s)}{\psi(z,s)} \,,
$$
where $\lambda,\mu \in \chr \tbT = \cochar \tbT^\vee$. 
Evidently, $\psi'$ is $q$-periodic in $z$ while 
\begin{equation}
\frac{\psi'(z,q^\sigma s)}{\psi'(z,s)}   = q^{-\lan \lambda, \sigma \ran
\lan \lambda, \mu \ran}  \, \mu(s)^{-\lan \lambda, \sigma \ran} \,
\lambda(s)^{-\lan \mu, \sigma \ran}\label{trans_1} \,. 
\end{equation}
Clearly, the  line bundle $\tau(\lambda \mu)^* \cU \otimes 
\cU^{-1}$ depends \emph{linearly} on $\lambda \mu$. This follows
either directly from \eqref{trans_1} or from basic general theorem about
translates of line bundles on abelian varieties.

\subsubsection{}

For ease of future reference, we note 
an immediate consequence of the formula \eqref{trans_1}.

\begin{Lemma}\label{lemma_UU} For $\lambda,\mu \in \chr \tbT = \cochar \tbT^\vee$ as
  above, the ratio 
$$
\frac{\vth(\lambda \cdot \mu)}{\vth(\lambda)\, \vth(\mu)} 
$$
is a meromorphic section of $\tau(\lambda \mu)^* \cU \otimes 
\cU^{-1}$\,. 
\end{Lemma}

\noindent 
Here, we had to distinguish between product $\lambda \cdot \mu$ of
functions on $\tbT$ and the composition
$$
\lambda \circ \mu \in \Hom(\tbT, \tbT^\vee) \,, 
$$
which enters the definition of $\tau(\lambda\mu)$.

We will be particularly interested in the case when 
$\mu$ is the weight $\hbar$ of the symplectic form. In this case, we
note the following 

\begin{Lemma}\label{lemma_UU2} The ratio 
$$
\prod \dfrac{\vth(\hbar w_i )}{\vth(w_i) \, \vth(\hbar)} 
$$
is a meromorphic section of 
$
\tau(\lambda \hbar )^* \cU \otimes 
\cU^{-1}
$ for  $\lambda = \prod w_i$. 
\end{Lemma}




\section{Elliptic stable envelopes}

\subsection{Attracting manifolds}

\subsubsection{} 

The setup is the same as e.g.\ Section 3.2 of \cite{MO1}. 
Let $\bA\subset \Ker \hbar \subset \bT$ be a subtorus. The 
normal weights to $X^\bA$ partition $\Lie \bA$ into finitely many
chambers. Let $\bAb \supset \bA$ be the toric compactification of 
$\bA$ defined by the fan of the chambers. A choice of a chamber 
$\fC$ defines a point 
\begin{equation}
0 = 0_\fC \in \bAb\label{0inbA}
\end{equation}
at infinity of $\bA$. For every $S\subset X^\bA$ we can define its 
attracting set 
$$
\Attr(S) = \{(x,s), s\in S, \lim_{a\to 0} a\cdot x = s \} \,,
$$
and the full attracting set $\fAttr(S)$ which is the minimal closed
subset of $X$ that contains $S$ and is closed under taking 
$\Attr(\,\cdot\,)$. Here and below a choice of the chamber $\fC$ 
is understood. 

\subsubsection{}\label{s_order_F} 

Let $F_i$ be the connected components of $X^\bA$.  We define
an ordering on the set of connected components by 
$$
F_1 \ge  F_2  \quad \Leftrightarrow \quad \fAttr(F_1) \cap F_2 \ne 
\varnothing  \,. 
$$
This is well defined because we assume the action to be linearized as
in \eqref{linearized}.

\subsection{Polarization}

\subsubsection{} 

A polarization of $X$, denoted 
$T^{1/2} X$, is an element of 
$K_\bT(X)$ such that 
\begin{equation}
TX =T^{1/2} X  + \hbar^{-1} \otimes 
\left(T^{1/2} X\right)^\vee \label{pol}
\end{equation}
where $\hbar$ is a character of $\bT$. We use the 
shorthand $T^{1/2}$ when $X$ is clear from context.

Polarization is a somewhat auxiliary but important technical concept that keeps
coming up in connection with stable envelopes. For instance,
in the elliptic
cohomology context, it is natural for elliptic
classes associated to Lagrangian
submanifolds of $X$ to be sections of line bundles closely related to
$\Theta(T^{1/2} X)$.

Cotangent bundles $X = T^* M$ 
have polarizations given by either base or the 
fiber direction, so that $\hbar$ is the weight of 
the natural symplectic structure on $T^* M$. Nakajima varieties come with natural polarizations \eqref{polNak}, 
which may be traced to their embedding in the cotangent 
bundle to a stack of quiver representations. 

With every polarization $T^{1/2}$ comes the opposite 
polarization $T^{1/2}_\textup{opp}$, which is the other term in 
the right-hand side of \eqref{pol}. 

\subsubsection{}

Suppose $\delta T^{1/2}$ is the difference of two polarizations in
$K_\bT(X)$. Then 
$$
\delta T^{1/2}  = \sum_{i=1}^l \left(w_i - \frac{1}{\hbar w_i} \right) 
$$
for certain Chern roots $w_i$, where each $w_i$ is a monomial in
equivariant variables and Chern roots of the tautological
bundles. Recall that $K_T(X)$ is a subscheme is the product of
$T$ with $S^{\bv_i} \C^\times$ and here we compute with coordinates
$w_i$ on
this ambient space. Denote 
$$
\lambda = \prod w_i  = \sqrt{ \hbar^{-l} \det \delta T^{1/2}} 
\in \Pic_T(X) \,. 
$$
By Lemma \ref{lemma_UU2}, we have
\begin{equation}
\Theta \left(\delta T^{1/2}\right) \cong 
\frac{\tau\left(-\lambda \hbar \right)^* \cU}{\cU} \otimes 
\Theta(\hbar)^{-l} \label{dT12} \,. 
\end{equation}
%

\subsubsection{}
The restriction of the polarization to $X^\bA$ can be decomposed 
\begin{equation}
T^{1/2}  \big|_{X^\bA}  = T^{1/2}  \big|_{X^\bA,>0} 
\oplus T^{1/2}  \big|_{X^\bA,\textup{fixed}} 
\oplus  T^{1/2}  \big|_{X^\bA,<0} \label{poltoA}
\end{equation}
into attracting, fixed, and repelling directions for the action of
$\bA$. 

The fixed part defines a polarization of $X^\bA$ which we denote 
by $T^{1/2} X^\bA$. We denote by
\begin{equation}
\ind =  T^{1/2}  \big|_{X^\bA,>0} \in K_\bT(X^\bA) \label{def_ind}
\end{equation}
the attracting part of polarization. This is an analog of the 
index of a fixed component. It is exactly the index if $X$ is a cotangent bundle
to some other manifold, the tangent bundle of which is given 
by $T^{1/2}$. 

\subsection{Definition of stable envelopes}

\subsubsection{} 
A $\bT$-equivariant map $f: Y \to X$ induces a diagram 
\begin{equation}
\xymatrix{
\bE_\bT(Y) \ar[d] &&  \Ell_\bT(Y) \times \cE_{\Pic_\bT(X)} \ar[d]
\ar[ll]_{1 \times f^*\quad} 
\ar[rr]^{\quad \Ell(f) \times 1}  && \bE_\bT(X) \ar[d]
\\
\cB_{\bT,Y} &&  \cB_{\bT,X} \ar[ll]_{1\times f^*} \ar[rr]^1 && \cB_{\bT,X} \,,
}\label{XYdiag}
\end{equation}
where $f^*$ is the pull-back of line bundles from $X$ to $Y$.

\subsubsection{} \label{s_defStab0} 

Stable envelopes is a way to map, with some shifts and twists, 
the universal bundles in \eqref{XYdiag}, as sheaves on the base,
for the inclusion of the fixed locus 
$$
\iota: X^\bA \to X \,, \quad \bA \subset \bT \,. 
$$
We have 
$$
\det \ind  \in \Pic_\bT(X^\bA) 
$$
and this defines a translation 
$$
\tau(-\hbar \det \ind) : \cB_{\bT,X^\bA} \to \cB_{\bT,X^\bA} \,. 
$$
We denote 
\begin{equation}
\cU' = (1 \times \iota*)^* \, \tau(-\hbar \det \ind)^* \,
\cU_{\bE_T(X^\bA)} \label{cUprime} \,. 
\end{equation}
This is a line bundle on the top left space in the following diagram
\begin{equation}
\xymatrix{
 \Ell_\bT(X^\bA) \times \cE_{\Pic_\bT(X)} \ar[d]
\ar[rr]^{\quad \Ell(\iota) \times 1}  && \bE_\bT(X) \ar[d]
\\
\cB_{\bT,X} \ar[rr]^1 && \cB_{\bT,X} \,. 
}\label{XYdiag2}
\end{equation}
Elliptic stable envelope is a map 
of $\cO_{\cB_{\bT,X}}$-modules 
\begin{equation}
\Theta\left(T^{1/2} X^\bA \right) \otimes 
\cU' 
\xrightarrow{\,\, \Stab_\fC \, \, } 
\Theta\left(T^{1/2} X\right) \otimes 
\cU \otimes \dots  \,,\label{Stab_map}
\end{equation}
where dots stand for meromorphic sections of a 
certain line bundle pulled back 
from 
$$
\cB' = \cB_{\bT,X} / \cE_\bA \,. 
$$
Equation \eqref{Stab_map} specifies the factors
of automorphy in the variables 
$a\in\bA$. The variables in $\cB'$, that is, the
rest of the equivariant and all K\"ahler variables, enter as 
parameters into this specification. The dependence on those 
variables is uniquely 
determined by a certain triangularity and normalization of 
stable envelopes, see Proposition \eqref{p_normalh} below. 

In the general spirit of stable 
envelopes, we constrain the $\bA$-dependence explicitly, and 
the rest is fixed by a certain uniqueness. In our personal
experience, this is the productive way to think about
stable envelopes and we ask those readers who feel
uncomfortable about the ellipsis in \eqref{Stab_map} to read on
until \eqref{Stab_map2}.

\subsubsection{}\label{s_def_support} 

Definition of stable envelopes involves supports, which means
the following. 
Let $s$ be a section of a coherent sheaf on $\Ell_\bT(X)$ over 
an open set in the base $\cB_{\bT,X}$ and let 
$$
f: Y\to X
$$
be an inclusion of a $\bT$-invariant set. We say 
$$
\supp(s) \subset Y  \quad \Leftrightarrow \quad f_{\textup{complement}}^*(s) = 0 
$$
where 
$$
f_{\textup{complement}}: \Ell_{\bT}(X\setminus Y) \to \Ell_{\bT}(X) 
$$
is the functorial map.  

\subsubsection{}\label{s_ind_support} 

In $X$, we have a decreasing sequence of closed sets 
\begin{equation}
  \label{defY_i}
   Y_i = \bigcup_{F_j \le  F_i} \Attr(F_j)\,, \quad Y_\infty=
   \varnothing\,, 
\end{equation}
and a typical strategy of proving that $s=0$ 
will be to show inductively that $\supp(s) \subset Y_i$ for all $i$. 

\subsubsection{}\label{s_defStab}

By definition, $\Stab_\fC$ satisfies the following two conditions:
\begin{itemize}
  \item[($\star$)] The support of $\Stab_\fC$ is triangular with respect to 
$\fC$, that is, 
 if locally over $\cB_{\bT,X}$  an elliptic cohomology 
class $s$ is supported on a component
$F_i$ of 
a fixed 
locus then $\Stab_\fC(s)$ is supported on $\fAttr(F_i)$. 
\item[($\star\star$)] Near the diagonal, we have 
  \begin{equation}
\Stab_\fC = (-1)^{\rk \ind} \, j_* \pi^* \,,\label{stab_diag}
\end{equation}
where 
$$
F_i  \xleftarrow{ \,\, \pi \, \,} \Attr(F_i) \xrightarrow { \,\, j
  \, \,} 
X
$$
are the natural projection and inclusion maps. Here, near 
the diagonal means that we restrict to the complement of $\bigcup_{F_j < F_i}
\Attr(F_j)$. 
\end{itemize}

\subsubsection{}

Note that the shift $\tau(-\hbar \det \ind)$ is necessary to make 
property ($\star\star$) agree with \eqref{Stab_map}. Indeed, let $N_{X^\bA,<0}$ 
be the repelling part of the normal bundle to $X^\bA$ so 
that 
$$
\Theta(-N_{X^\bA,<0})  \xrightarrow{\quad j_* \pi^* \quad} 
\cO_{\Ell_\bT(X)} \,.
$$
We observe from Lemma \ref{lemma_UU2} that 
$$
\Theta\left(T^{1/2} X - T^{1/2} X^\bA - N_{X^\bA,<0}+
(\rk \ind ) \, \hbar \right)
\Big|_{\textup{$\cE_\bA$-orbits}} \cong 
\frac{\tau(-\hbar \det \ind)^* \cU}{\cU} \,, 
$$
or, in other words, that the ratio of normal weights 
$$
\prod_{w\in \textup{weights}(\ind)}
\frac{\vth(w) \, \vth{(\hbar)}}{\vth(\hbar^{-1} w^{-1})}
$$
is a section of $\cU^{-1} \otimes \tau(-\hbar \det \ind)^* \cU$. 
Here $\rk \ind$ is the rank of \eqref{def_ind}.

\subsubsection{}
Stable envelopes depend on the choice of polarization, but 
as \eqref{dT12} shows stable envelopes for different polarizations
are related by a shift of K\"ahler parameters. 

Also note that one of the factors of the 
base $\cB_{\bT,X}$ from \eqref{bBT} is naturally an extension
$$
0 \to \cE_\bT^\vee \to \cE_{\Pic_\bT(X)} \to \cE_{\Pic(X)} \to 0  \,. 
$$
As a simple corollary of the uniqueness of stable envelopes, see
below, they are constant sections of trivial bundles 
along the $\cE_\bT^\vee$-orbits on the base $\cB_{\bT,X}$. 
In particular, we don't introduce any special symbols to denote
variables in $\cE_\bT^\vee$ since nothing depends on them. 

It is, however, convenient to keep these directions because the shifts 
\eqref{def_tau} involve $\lambda\in \Pic_\bT(X)$ and really change 
if $\lambda$ is twisted by a character of $\bT$. 
A concrete form of this dependence is given by 
Lemma \ref{lemma_UU}.


\subsubsection{}

In Theorem \ref{p_uni} in Section \ref{s_unique} below we will show
that the properties required in Sections \ref{s_defStab0} and 
\ref{s_defStab} determine stable envelopes uniquely. 
In particular, the dots in \eqref{Stab_map} mean certain 
automorphy factors in all variables other than those in
$\cE_\bA$. By uniqueness, these can be read off from the 
restriction 
\eqref{stab_diag} of stable envelopes to the diagonal. 

\begin{Proposition}\label{p_normalh} Assuming uniqueness, elliptic stable 
envelopes give a map 
\begin{equation}
\Theta\left(T^{1/2} X^\bA \right) \otimes 
\cU' \otimes \Theta(\hbar)^{-\rk \ind} 
\xrightarrow{\,\, \Stab_\fC \, \, } 
\Theta\left(T^{1/2} X\right) \otimes 
\cU \,,\label{Stab_map2}
\end{equation}
where the shift $\cU'$ is a shift of the universal bundle as  in
\eqref{cUprime} and $\rk \ind$ denotes the rank of \eqref{def_ind}. 
\end{Proposition}

\begin{proof}
  Follows from Lemma \ref{lemma_UU2}. 
\end{proof}

\subsection{Example: $T^*\bP(W)$}\label{seTP} 

\subsubsection{} 
Let $W$ be a vector space and consider the 
$GL(W)\times GL(1)$ module 
\begin{equation}
M = W \otimes \C \,, \label{eqMW}
\end{equation}
where $\C$ denotes the defining representation of $GL(1)$. 
We have 
\begin{align*}
X &= T^*\bP(W)\\
 & = \{(v,\xi)\in \rM \oplus \rM^*, \lan \xi,v \ran =0, v\ne 0 \} /GL(1) \,.
\end{align*}
where
$$
s \cdot (v,\xi)= (s v, s^{-1} \xi) \,, \quad s \in GL(1) \,. 
$$
In addition to $GL(W)$, the group $\C^\times_\hbar$ acts on $ \rM \oplus \rM^*$ and 
$X$ by 
$$
\hbar \cdot (v,\xi) = (v,\hbar^{-1} \xi) \,, \quad \hbar \in
\C^\times_\hbar\,. 
$$
We denote $G=GL(W) \times \Ct_\hbar$. 

\subsubsection{}\label{s_cM} 
The vector space $\rM$ descends to a bundle 
$\cM\cong \cO(1)^n$ on $X$ and we take 
\begin{equation}
T^{1/2} X = \cM - \cO_X  \label{polPn}
\end{equation}
which is the pullback of the tangent bundle 
under $X \to \bP(W)$. 

\subsubsection{} 
Let 
$\bT\subset G $ be a maximal torus and $\bA = \bT \, \cap \, GL(W)$. 
The torus $\bA$ acts diagonally in some basis $\{e_i\}\subset W$ with 
weights $a_i$. Its fixed points are 
$$
F_i = \{ \C e_i\} \in \bP(W) \subset X \,,
$$
and we assume that $\fC$ and the labeling of $\{F_i\}$ are such that 
$$
F_1 > F_2 > \dots > F_n \,.
$$

\subsubsection{} 
The line bundle $\cO(1)$ associated to the fundamental 
weight $s$ of $GL(1)$ generates $K_\bT(X)$ and has weight $a_i^{-1}$ 
at the $i$th fixed point. Hence 
$$
\spec K_\bT(X) = \left\{ \prod(1-s a_i) =0 \right\} \subset \bT \times GL(1) \,.
$$
The reduction of this modulo $q^\Z$ is $\Ell_\bT(X)$.  We continue
to use the product 
to denote the group operation on the quotient. 

The coordinate $s$ on $GL(1)/q^\Z = E=E^\vee$ is the elliptic $c_1(\cO(1))$ and we
denote by $z\in E$ the dual K\"ahler parameter.

\subsubsection{}\label{s_st_Pn}
We will now check that the function  
\begin{equation}
\Stab(F_k)  = \prod_{i<k} \vth(sa_i) \,
\frac{\vth(sa_kz\hbar^{k-n})}{\vth(z\hbar^{k-n})} \,
\prod_{i>k} \vth(sa_i\hbar)\label{StabPn}
\end{equation}
satisfies the definition for stable envelopes. 
Note a match between the terms in numerator and 
denominator of this expression to the terms in 
the polarization \eqref{polPn}. 

The first thing to check is that \eqref{StabPn} is a section of 
a correct line bundle when restricted to $\cE_\bA$-orbits in the
base. This means that as function of $a_i$'s and $s$ it has
to have the same factors of automorphy as the following 
product: 
\begin{alignat}{2}
  &\prod_{i=1}^n \vth(sa_i)\times&&\textup{from polarization} \label{prodprod}\\
  & \frac{\vth(sz)}{\vth(s) \vth(z)}\times &&\textup{for $\cU$ in the
    target}\notag\\
& \frac{\vth(a_k^{-1}) \vth(z)}{\vth(a_k^{-1} z)}\times &&\textup{for $\cU$ in the
    source}\notag\\
  & 
  \frac {\vth( a_k^{k-n} \prod a_i) \vth(\hbar)}
  {\vth(\hbar^{-1} a_k^{k-n} \prod a_i)}  &\qquad& 
\textup{for $\tau(-\hbar\det\ind)$ in the source} \,, 
\notag 
\end{alignat}
which is indeed the case. The factors in \eqref{prodprod} have 
the following explanation. 

Polarization in the target of $\Stab$ is given by \eqref{polPn} 
while the polarization in the source is trivial, thus the first line 
in \eqref{prodprod}. The variable $s$ is defined as the dual to 
the coordinate $z$ on $\cE_{\Pic(X)}$, which, in turn, corresponds
to $\cO(1)$ with the chosen  linearization. Thus the second line in 
\eqref{prodprod} is obtained by copying \eqref{cUthth}. 
Restricted to $F_k$, $\cO(1)$ becomes the trivial bundle with 
weight $a_k^{-1}$, thus the third line in \eqref{prodprod} is the 
reciprocal of the second with the substitution $s=a_k^{-1}$. Finally 
$$
\ind(F_k) = \sum_{i>k}  \frac{a_i}{a_k} 
$$
therefore 
$$
\det\ind(F_k) = a_k^{k-n} \prod_{i>k} a_i 
$$
which gives the fourth line in \eqref{prodprod} as in
Lemma \eqref{lemma_UU2}.

\subsubsection{}
Evidently, the function \eqref{StabPn} vanishes when 
restricted to $F_i$ with $i<k$ because of the factors
$\theta(sa_i)$. When we restrict to $F_k$, we get 
precisely the product of $\theta$-functions of the 
repelling weights, with the correct sign.

\subsection{Uniqueness}\label{s_unique}

\begin{Theorem}\label{p_uni} 
Elliptic stable
envelopes are unique. 
\end{Theorem}

The key logical point in the argument below will be the following rigidity
statement. If $s\ne 0$ is a regular section of a degree zero line
bundle $\cL$ on
an abelian variety $\cE$ then $\cL=\cO_\cE$ is the trivial line
bundle. For us, $\cE$ will be an $\cE_\bA$-orbit in $\cB_{X,\bT}$ and
$\cL$ a certain relative of the universal bundle $\cU$ restricted to
$\cE$. The latter has a nontrivial $z$-dependence, and in particular
$\cL \not\cong \cO_\cE$. This rules out nonzero sections of $\cL$. 

\begin{proof}
Suppose for some $\fC$ there are two different maps $\Stab$ and 
$\Stab'$ satisfying the above conditions. Let $F_i$ be a component 
of the fixed locus and consider the map 
$$
\delta = \big(\Stab-\Stab'\big) \Big|_{\textup{$F_i$ summand in 
the source}} \,. 
$$
By the normalization condition
$$
\supp(\delta) \subset \bigcup_{F_j < F_i} \Attr(F_j)  \,,
$$
and we will argue as in Section \ref{s_ind_support} that, in fact, 
$\supp(\delta)=\varnothing$. 

The closed sets $Y_j$ in \eqref{defY_i}
form a partially ordered set and we fix  a maximal element $Y_k$ 
among those that intersect $\supp(\delta)$. We set 
\begin{align*}
  X' &= X \setminus \bigcup_{F_j \not\ge F_k} \Attr(F_j)  \\
  Y' &= X' \cap Y_k  \,. 
\end{align*}
{}The long exact sequence of the pair $(X',X'\setminus Y')$
 in elliptic cohomology starts out as 
\begin{equation}
0 \to \Theta(-N_{X'/Y'}) \to  \cO_{\Ell_\bT(X')} \to 
\cO_{\Ell_\bT(X'\setminus Y')} \to \dots  \label{longell} 
\end{equation}
and by our hypothesis $\delta$ restricts to zero on $X'\setminus Y'$. 
Our goal is to show that $\delta$ restricts to zero on $X'$. 
Since $Y'$ is a vector bundle over the $F_k$, it is enough 
to show that the pullback of $\delta$ under 
$$
\iota_k : F_k \to X
$$
vanishes. 

Let 
$$
s = \iota^*_k \,\delta \Big|_{\cE_{\bA} b}
$$
be the restriction of this pullback to a general $\cE_\bA$-orbit
$$
\cE_{\bA} b  \subset \cB_{\bT,X}\,, \quad b=(t,z) \in \cE_\bT \times 
\cE_{\Pic_\bT(X)}
$$
in the base $\cB_{\bT,X}$. Since $\bA$ does not act on 
$F_k$, the restriction of $\cO_{\Ell_\bT(F_k)}$
to $\cE_{\bA} b$ is a trivial bundle 
with fiber $\Hd(F_k^{\bT_t},\C)$, where
$\bT_t$ is as Section \ref{s_flatness}.
 Therefore, $s$ is a regular section of a
trivial bundle twisted by a line bundle. 

{}From \eqref{longell}, we see that $s$ is divisible by 
$$
\theta\in H^0\left(\Theta\big(TX\big|_{F_k,<0}\big)\right)\,, 
$$
where $TX\big|_{F_k,<0}$ are the repelling directions in the normal
bundle to $F_k$. The $\bA$-weights in $TX\big|_{F_k,<0}$ are, 
up to a sign, the same as the $\bA$-weights of $T^{1/2}X$. 
Therefore
$$
\theta^{-1}  s 
\in H^0( \cE_\bA \, b, \cO^{\rk} \otimes \cL_{t,z})\,, \quad \deg \cL_{t,z}
= 0 \,,
$$
where $\cL_{t,z}$ is a certain combination of Thom bundles and 
the universal bundle. Since it is a line bundle of degree zero 
$$
s \ne 0 \quad \Longrightarrow \quad \textup{$\cL_{t,z}$ is trivial} \,. 
$$
This is, however, impossible
because $\cU\big|_{\cE_\bA}$, and hence $\cL_{t,z}$, 
 depends 
nontrivially on the variable $z$. Indeed, since $F_i > F_k$, there 
exists a chain of $\bA$-invariant rational curves that flows from $F_k$ to 
$F_i$. An ample line bundle on $X$ will have a positive degree on this 
chain of curves, hence different $\bA$-weights on $F_i$ and $F_k$. 
Since there is a line bundle with different weights on $F_i$ and
$F_k$, the dependence on $z$ is nontrivial. So, 
$$
\textup{$\cL_{t,z}$ is nontrivial} \quad \Longrightarrow \quad 
s =0 \quad \Longrightarrow \quad \iota^*_k \,\delta =0 
\quad \Longrightarrow \quad \delta =0  \,.
$$
\end{proof}

\begin{Corollary}\label{c_uni} Stable envelopes preserve support in
  the sense that 
$$
\supp \Stab(\,\cdot\, ) \subset \fAttr \left( \supp(\,\cdot\,)\right)
\,. 
$$
\end{Corollary}

\begin{proof}
Apply the logic above to $X \setminus \fAttr \left(
  \supp(\,\cdot\,)\right)$. 
\end{proof}

\subsection{Triangle lemma}

\subsubsection{}\label{sfC'}

Let $\fC'\subset \fC$ be a face of some dimension, as in Section 3.6 
in \cite{MO1} and let $\bA' \subset \bA$ be the subtorus associated 
to the span of $\fC'$ in $\Lie A$. We have a triangle of 
embeddings 
$$
\xymatrix{
X^\bA \ar[rr] \ar[rd]&&X \\
&X^{\bA'} \ar[ru] 
} 
$$
to each arrow in which we can associate a stable envelope map. 
In particular, we have the stable envelope map $\Stab_{\fC/\fC'}$, 
which we can pull back to $\cB_{\bT,X}$ using the composition 
of the maps 
$$
\boldsymbol{\tau} : \cB_{\bT,X}  
\to 
\cB_{\bT,X^{\bA'}} 
\xrightarrow{\,\tau(-\hbar\det \ind)\,} 
\cB_{\bT,X^\bA} \to 
\cB_{\bT/\bA,X^\bA}
$$
where the first and the last are the natural maps and the middle 
translates by 
$$
\det \ind \in \Pic_\bT(X^{\bA'})\,.
$$
Here the index is the index associated to the embedding $X^{\bA'} \to
X$ by the given polarization of $X$. We denote
$$
\Stab_{\fC/\fC'}(\, \cdot \, - \hbar\det \ind_{X^{\bA'}}) = 
\boldsymbol{\tau}^* \Stab_{\fC/\fC'} \,.
$$

\begin{Proposition}\label{p_tria} 
We have
$$
\Stab_\fC  = \Stab_{\fC'} \circ 
\Stab_{\fC/\fC'}(\, \cdot \, - \hbar\det \ind_{X^{\bA'}})\,. 
$$
\end{Proposition}

\begin{proof}
Immediate from the uniqueness of stable envelopes. 
\end{proof}

\subsubsection{Example} 
For $T^*\bP(W)$, see section \ref{seTP} we can 
take
$$
W' = \bigoplus_{i=1}^m \C e_i \subset W
$$
and $\bA' \subset \bA$ the point-wise stabilizer of $\bP(W') 
\subset \bP(W)$, so that 
$$
X^{\bA'} = T^* \bP(W')  \sqcup \{F_{m+1},\dots,F_n\} \,. 
$$
Then 
$$
\ind_{T^* \bP(W')} \cong \sum_{i={m+1}}^n a_i \cO(1) 
$$
and hence the shift of the stable envelope by $-\hbar \det \ind$ 
is the shift of the K\"ahler parameters by $-\hbar \cO(n-m)$, 
which means $z\mapsto z \hbar^{m-n}$. These are precisely the 
shifts in the formula \eqref{StabPn}.

\subsection{Duality}\label{s_dual}

\subsubsection{}

Let $f: X \to Y$ be a $\bT$-equivariant map between formal
$\bT$-varieties and suppose that the restriction of $f$ to 
$X^\bT$ is proper. Equivariant localization lets one define $f_*$ 
as the unique map completing the diagram 
\begin{equation}
\xymatrix{
\Theta_{X^\bT}\left(T X^\bT-f^*TY\right)  \ar[rr]^{\iota_*} 
\ar[dr]_{(f\,\circ \,\iota)_*}&& 
\Theta_{X}\left(T X - f^*TY\right) \ar[dl]^{f_*}\\
&\cO_{\Ell_\bT(Y)} 
} \,\, .\label{eqf*}
\end{equation}
In particular, for the map $f: X\to \pt$ we get 
 the composite map 
\begin{equation}
 \Theta\big(T^{1/2}X\big) \otimes \Theta\big(T^{1/2}_\textup{opp}
  X\big) 
\xrightarrow{\textup{ product }} \Theta(TX) 
 \xrightarrow{\,\, f_* \,\,}\cO_{\cE_\bT,\textup{localized}} \,,
\label{dualTh} 
 \end{equation}
in which the first map is multiplication and the second 
pushforward.

\subsubsection{}

Define
$$
\Stab^*_{\fC,T^{1/2}X} (z) = 
\Stab_{-\fC,T^{1/2}_\textup{opp} X} (-z)^\vee
$$
where $z\in \cE_{\Pic_\bT(X)}$ is the K\"ahler parameter and 
the dual means transpose with respect to dualities 
\eqref{dualU} and \eqref{dualTh}. This is a meromorphic map 
\begin{equation}
\Theta\left(T^{1/2} X\right) \otimes 
\cU
\xrightarrow{\,\, \Stab^* \, \, } 
\Theta\left(T^{1/2} X^\bA \right) \otimes 
\cU' 
\otimes \dots  \,,\label{Stab_map3}
\end{equation}
where, again,  dots stand for a 
certain line bundle pulled back 
from 
$\cB' = \cB_{\bT,X} / \cE_\bA$ \,. 

\subsubsection{}
Duality for stable envelopes is the following 

\begin{Proposition}\label{p_dual}
$$
\Stab^* \circ \Stab = 1  \,. 
$$
\end{Proposition}

\begin{proof}
Since the attracting and repelling manifolds intersect properly, 
the composition is regular. 

Fix two components $F_i$ and $F_j$ and restrict $\Stab^* \circ
\Stab $ to the corresponding direct summands in the source and 
in the target. By the support condition, this restriction is only nonzero if $F_i \ge
F_j$.

If $F_i > F_j$ then the argument already used in the proof 
of Theorem \ref{p_uni} shows the universal bundle $\cU$ has different 
$z$-dependence in the source and in the target. It follows that 
$\Stab^* \circ \Stab$ vanishes in this case also. 

For $F_i=F_j$, condition $(\star\star)$ in Section \ref{s_defStab} 
implies $\Stab^* \circ \Stab $ is the identity, which concludes the
proof. 
\end{proof}

\section{Existence of stable envelopes}
\label{s_exi}

\subsection{Hypertoric varieties}

\subsubsection{}

Hypertoric varieties $X$, see e.g.\ \cite{Proud}, 
are algebraic symplectic reductions of a vector 
space by an action of a torus $S$, that is, 
$$
X = T^*\rM \rd  S  = \mu^{-1}(0) \rdd S = \mu^{-1}(0)_\textup{ss} / S 
$$
where $\rM$ is a representation of $S$ and 
$$
\mu: T^* \rM \to \Lie(S)^* 
$$
is the (holomorphic) moment map. As the torus $\bT$ we can take 
a maximal torus 
$$
\bT \subset  \C^\times_\hbar \times GL(\rM)^S\,,
$$
where $\C^\times_\hbar$ acts by rescaling the cotangent directions 
with weight $\hbar^{-1}$.

\subsubsection{}

An example of a hypertoric variety is  $T^*\bP(\rM)$ discussed
 in Section \ref{seTP}. In this case $S\cong \C^\times$ acts on $\rM$
by scalars. 

In this paper, we assume the $S$-action on $\mu^{-1}(0)_\textup{ss}$ is free and so 
$X$ is smooth.  This is a very strong restriction on the 
$S$-module $\rM$. Namely, if we think of weights of $\rM$ as a matrix 
$$
\Z^{\dim \rM} \to \chr(S) \cong \Z^{\rk S}
$$
then this matrix is surjective and determinants of 
all its submatrices (in particular, all matrix elements) are in 
$\{0,\pm 1\}$. 

\subsubsection{}\label{s_hypertoric_prod} 

Generalizing the explicit formula \eqref{StabPn}, we have the 
following: 

\begin{Proposition}
Stable envelopes exist for hypertoric varieties.
\end{Proposition}

\begin{proof}
We first deal with the case $\bA = \Ker \hbar$, in which case 
$X^\bA$ is a finite set. Let $F$ be one of the fixed points. 
We may assume 
$$
F = T^*\rM_0 \rd S
$$
where the decomposition 
\begin{equation}
\rM = \rM_0 \oplus \rM_1 \label{M01}
\end{equation}
is such that 
\begin{itemize}
\item $S$ acts as a maximal torus on $\rM_0$ 
$$
S\cong \left \{ 
  \begin{pmatrix}
    s_1  & \\
     & s_2 \\
         &  & \ddots  
 \end{pmatrix} 
\right\} \subset GL(\rM_0) \,, 
$$
\item 
$\bA$ acts as a maximal torus in $\rM_1$ and with weights 
$\alpha_i$ on $\rM_0$ \,. 
\end{itemize}
Further we may assume that 
$$
T^{1/2} X = \cM - \Lie(S)
$$
where $\cM$ is the bundle associated to the representation $\rM$ and 
$\Lie(S)$ is the trivial bundle with this fiber. In 
accordance with \eqref{M01}, $\cM=\cM_0\oplus 
\cM_1$ and we can further assume
that 
$$
T^{1/2}\big|_{F,<0} = \cM_{1}  \,.
$$
Then 
\begin{equation}
\Stab(F) = \vth(\cM_1) \, \prod_{i=1}^{\rk S} \, \frac{\vth(s_i z_i
  \alpha_i)}{\vth(z_i)}\label{s_theta_prod}
\end{equation}
where $\vth(\cM_1)$ is the product $\vth(x_i)$ over the Chern roots
$x_i$ of $\cM_1$. This generalizes formula \eqref{StabPn}. 

For subtori $\bA' \subset \bA$, we can use Proposition 
\ref{p_tria} and duality from Proposition \ref{p_dual} to 
write 
$$
\Stab_{\fC'} = \Stab_\fC  \circ \Stab_{\fC/\fC'}(\, \cdot \, - \hbar\det
\ind_{X^{\bA'}})^* \,. 
$$
The left-hand side of this formula defines stable envelopes for 
general subtori in $\Ker \hbar$. 
\end{proof}

\subsection{Nakajima varieties}

\subsubsection{} 

Now suppose 
\begin{equation}
X = T^*\rM \rd  G  = \mu^{-1}(0) \rdd G \label{XrdM}
\end{equation}
where $G \subset GL(\rM)$ is reductive and 
$$
\mu: T^*\rM \to \fg^*\,, \quad \fg = \Lie G \,, 
$$
is the moment map for the $G$-action. If we additionally assume that
the 
$G$-action on $\mu^{-1}(0)_\textup{ss}$ is free and $G$ is nonabelian, 
then the only examples of such quotients known to us 
are Nakajima varieties. 

\subsubsection{}\label{s_polsG}

As before, we take $\bT$ to be the maximal torus 
$$
\bT \subset \C^\times_\hbar \times GL(\rM)^G \,,
$$
and as a polarization we can take 
\begin{equation}
T^{1/2} X = \cM - \cG  \,,\label{polNak}
\end{equation}
where $\cG$ is the bundle on $X$ associated to the adjoint 
action of $G$ on $\fg$. 

\subsubsection{}

\begin{Theorem}\label{t1} 
Elliptic stable envelopes exist for Nakajima varieties. 
\end{Theorem}

The proof will reduce the statement to the existence of stable 
envelopes for hypertoric varieties. The general circle of ideas relating 
the cohomology of a $G$-quotient to the cohomology of the 
quotient by the maximal torus is known as \emph{abelianization}, 
see in particular \cite{HausProud}.  Abelianization of stable envelopes in 
the usual equivariant cohomology was developed by 
D.~Shenfeld \cite{Shen}, see also the exposition in 
\cite{Smir1}. 

\subsection{Proof of Theorem \ref{t1}} 

\subsubsection{} 

Let $S$ be a maximal torus of $G$ and let $B\supset S$ be a Borel 
subgroup with Lie algebra $\fb$. 
Choose a maximal compact subgroup $U\subset G$ and a $U$-invariant 
Hermitian metric $\|\, \cdot\, \|$ on $\rM$. This defines a real moment 
map 
$$
\mu_\R: T^*\rM \to \fu^*\,, \quad \fu= \Lie U  \,. 
$$
Consider the diagram 
\begin{equation}
  \label{diagF}
  \xymatrix{
\fF \ar[d]_\pi \ar[rr]^{\bj_+\quad\quad}&&
 \mu^{-1}(\fb^\perp)\rdd S \ar[rr]^{\quad\quad \bj_-}&& X_S \\
X
}
\end{equation}
in which 
\begin{equation}
  \label{eq:8}
  \fF = \mu_\R^{-1}(\eta)/ \left( U \cap S\right)\,, \quad 
  X_S = T^*\rM \rd  S \,,
\end{equation}
where the parameter $\eta \in (\fu^*)^U$ corresponds to 
the choice of the stability condition for the GIT 
quotient\footnote{Nakajima denotes the parameter $\eta$ by 
$\theta$.}.  

The fibers of $\pi$ are flag varieties
\begin{equation}
  \label{fl_iso}
 U/(U \cap S) \cong G/B \,. 
\end{equation}
The choice of $B$ in this isomorphism will be important below.

\subsubsection{}\label{sfF} 

While $\fF$ is not a complex subvariety of $X_S$, the map 
$\pi$ is a $G/B$-bundle and therefore is complex-oriented. 
The description of elliptic cohomology of $G/B$-bundles, 
see Section 1.9 in \cite{GKV}, 
shows: 
\begin{itemize}
\item The cohomology of $\fF$ is generated as a module over 
the cohomology of the base by tautological line bundles. In 
particular, it is tautological if  the cohomology of $X$ is
tautological. 
\item The map of sheaves
$$
\pi_* : \Theta(\Ker d\pi) \to \cO_{\Ell_\bT(X)} 
$$
is surjective. 
\end{itemize}

\subsubsection{} 
Recall the polarization of $X$ introduced in 
Section \ref{s_polsG}. 
We have
$$
\pi^* \cG = \cN \oplus \cN^\vee \oplus \Lie(S) \,,
$$
where $\cN$ is the bundle on $X_S$ 
associated to the adjoint action of $S$ on $\fn = [\fb,\fb]$
and $\Lie(S)$ is a trivial bundle with this fiber. 

The isomorphism \eqref{fl_iso} gives
\begin{equation}
\Ker d\pi \cong \cN^\vee \,. \label{kerdpi}
\end{equation}
Therefore we may replace $\pi^* \cG$ by 
$\Ker d\pi \oplus \hbar^{-1} \cN^\vee \oplus \Lie(S)$ in a 
polarization. This gives a polarization of $T^{1/2} X_S$ such 
that 
\begin{equation}
\bj^*\left(T^{1/2} X_S - \hbar^{-1} \cN^\vee\right) = \pi^*(T^{1/2} X)
+ \Ker d\pi \,,\label{polpull}
\end{equation}
where $\bj = \bj_- \circ \bj_+$. 

\subsubsection{} 
Consider the following chain of maps 
\begin{multline}
  \label{chn}
\Theta(T^{1/2} X) \xleftarrow{\quad \pi_* \quad} 
\Theta(\pi^*(T^{1/2} X)+\Ker d\pi) \xleftarrow{\quad\,\, \bj_+^*
  \quad}  \\
\xleftarrow{\quad\quad} \Theta(\bj_-^*(T^{1/2}
X_S-\hbar^{-1} \cN^\vee))
\xrightarrow{\quad\,\, {\bj_{-,*}}\quad} 
\Theta(T^{1/2} X_S) \,,
  \end{multline}
in which we used the identification 
\begin{equation}
\bj_-^* \, \hbar^{-1} \cN^\vee = \textup{normal bundle to $\bj_-$}
\,. \label{normjm}
\end{equation}
We claim the map $\bj_+^*$ in \eqref{chn} is also surjective. 
Indeed any tautological class extends to a class on all of $X_S$. 

\subsubsection{} 
The map 
$$
\Pic(X) = \chr G \to \chr S = \Pic(X_S) 
$$
induces an embedding $\cB_{\bT,X} \hookrightarrow \cB_{\bT,X_S}$ 
such that all maps in \eqref{chn} are defined after tensoring 
with the universal bundle and restricting to the image of this 
embedding. 

\subsubsection{}

Let $\bA$ be a torus in the kernel of $\hbar$ and let $F$ be a 
component of $X^\bA$. The torus $\bA$ acts on 
$\pi^{-1}(F)$.  Let $F'$ be a component of $\pi^{-1}(F)^\bA$. 
The composite map $\pi'_*$ in the diagram 
\begin{equation}
  \label{pi'}
  \xymatrix{
F' \ar[r] \ar[dr]_{\pi'} & \pi^{-1}(F) \ar[d]^\pi \\
& F 
}
\end{equation}
is still surjective on elliptic cohomology because it is 
still a fibration by flag varieties.  

\subsubsection{}\label{s_choiceF} 

Of all possible $F'$ we pick the one for which the 
normal weights to $F'$ in 
$\pi^{-1}(F)$ are \emph{repelling} for the action of $\bA$, 
or in other words 
\begin{equation}
  \label{Ker>0}
  \left(\Ker d\pi \big|_{F'}\right)_{>0} = 0 \,. 
\end{equation}
Recall that the fibers of $\pi$ are smooth and connected, therefore
there is exactly one component of the fixed locus for which all normal
weights are repelling. 

\subsubsection{}
The points of $F'$ are fixed by $\bA$ on the quotient by $S$, which 
means that there is a map 
$$
\phi_F: \bA \to S 
$$
so that the preimage of $F'$ in $T^*\rM$ is fixed point-wise under 
the action of $(a,\phi(a)) \in \bA \times S$. This induces an 
action of $\bA$ on $\fg$ and \eqref{Ker>0} means 
\begin{equation}
\left(\fn^\vee\right)_{>0} = 0 \label{fn>0}
\end{equation}
or equivalently 
\begin{equation}
  \label{N>0}
  \left(\cN^\vee \big|_{F'}\right)_{>0} = 0 \,. 
\end{equation}

\subsubsection{} 

Consider the following analog of \eqref{diagF} 
\begin{equation}
  \label{diagF2}
  \xymatrix{
F' \ar[d]_{\pi'} \ar[rr]^{\bj'_+\quad\quad}&&
F_S \cap \mu^{-1}(\fb^\perp)\rdd S \ar[rr]^{\quad\quad \bj'_-}&& F_S \\
F
}
\end{equation}
where $F_S$ is the component of $X_S^\bA$ that contains $F'$. 
As before, $\left(\bj'_+\right)^*$ is surjective and 
the equations \eqref{polpull}, \eqref{Ker>0}, and \eqref{N>0} 
imply that 
\begin{equation}
\pi^*(\ind F)  = (\bj')^*\left(\ind F_S \right) \,,  \label{indpull}
\end{equation}
where $\bj' = \bj'_- \circ \bj'_+$.

\subsubsection{}

We want to define the stable envelope as the 
composition 
\begin{equation}
  \label{diagXs}
  \xymatrix{
\Theta\left(T^{1/2} F \right) \otimes 
\cU'  
\ar[rrr]^{\bj'_{-,*} \left((\bj'_+)^*\right)^{-1} \left(\pi'_*\right)^{-1}}  
 \ar[d]_{\Stab}
&&&
\Theta\left(T^{1/2} F_S\right) \otimes 
\cU' \ar[d]_{\Stab_S} \\
\Theta\left(T^{1/2} X \right) \otimes 
\cU 
&&& \Theta\left(T^{1/2} X_S \right) \otimes 
\cU 
\ar[lll]_{\pi_*\,  \bj_+^*\, (\bj_{-,*})^{-1}}
}
\end{equation}
in which the index shifts agree by \eqref{indpull}. 

We need to check that the inverse maps that 
appear in both horizontal arrows on 
\eqref{diagXs} are well defined. 

\subsubsection{}
The inverse in the top horizontal line of \eqref{diagXs} means 
that we pick a tautological 
class on $X_S$ that maps to the right class on $X$ under 
$\pi'_*  \circ (\bj'_+)^*$. While the corresponding map of 
sheaves is surjective, taking the preimage is by itself not a
well-defined operation. However, $\bA$ acts trivially on $F$ 
and therefore the corresponding sheaves are \emph{trivial} along the 
$\cE_\bA$-orbits on the base $\cB_{X,\bT}$. Therefore, the preimage may be 
chosen locally on $\cB_{X,\bT}/\cE_\bA$. The resulting stable
envelope maps glue by the uniqueness proven in 
Theorem \ref{p_uni}. 

\subsubsection{}

It remains to check that the inverse in 
$$
 (\bj_{-,*})^{-1} \circ \Stab_S \circ \, \bj'_{-,*}
$$
is well-defined. 

Since $\mu$ is an equivariant map and the normal 
weights to $\fb^\perp \subset \fg^*$ are non-attracting
by \eqref{fn>0} it follows that 
$$
\fAttr (\mu^{-1}(\fb^\perp)^\bA) \subset \mu^{-1}(\fb^\perp) \,. 
$$
{}From Corollary \ref{c_uni} we now conclude 
that $\Stab_S \circ \, \bj'_{-,*}$ factors through $\bj_{-,*}$
and hence the diagram \eqref{diagXs} is well-defined.

\subsubsection{} 

By construction, we get a section of the correct line bundle with the
correct support as in the condition
$(\star)$ in Section \ref{s_defStab}. To finish
the proof, we need to verify condition $(\star\star)$ in the same
section.

Consider the maps in the diagram \eqref{diagXs} in the neighborhood of
$F'$. The hypertoric map $\Stab_S$ satisfies $(\star\star)$ and,
therefore, we need to analyze the difference between the vector bundles
$\pi^*\left(N_{X/F}\right)_{<0}$ and $\left(N_{X_S/F'}\right)_{<0}$.
By
our analysis, new repelling normal directions in $X_S$
appear either in \eqref{kerdpi} or in \eqref{normjm}. In 
the first case, they cancel out upon $\pi_*$. In second case, they
cancel out upon $(\bj_{-,*})^{-1} \circ \Stab_S \circ \bj'_{-,*}$.
This concludes the proof.

\subsection{Example: $T^*\Gr(k,n)$} 

\subsubsection{}

The cotangent bundle $X=T^*\Gr(k,n)$ of the Grassmannian 
is obtained for 
$$
\rM = \C^n \otimes \C^k 
$$
with $G=GL(k)$ and 
$$
\bA  = 
\begin{pmatrix}
 a_1 \\
& \ddots \\
& & a_n 
\end{pmatrix}
\subset GL(n) \,. 
$$
The abelianization of $X$ 
is the product 
$$
X_S = T^* \bP(\C^n)^{\times k} 
$$
of cotangent bundles of projective spaces discussed in 
Section \ref{seTP}. We choose the chamber $\fC$ as in 
that example, which means that 
$$
a_j / a_i > 0   \quad \Leftrightarrow \quad  i < j \, .
$$

\subsubsection{}
A particular simplifying feature of this example is that
the $\bA$-fixed points on $X_S$ are isolated. 
Concretely, the points in $X^\bA$ correspond to coordinate
subspaces in $\Gr(k,n)$ and are indexed by $k$-element 
subsets of $\{1, \dots, n\}$. The fixed points above
them correspond to injective maps 
\begin{equation}
\mu : \{1,\dots,k\} \to \{1, \dots, n\} \,. \label{mapmu}
\end{equation}
Other points in $X_S^\bA$ are $G$-unstable.

\subsubsection{}
We identify coordinates in 
$$
S = \begin{pmatrix}
 s_1 \\
& \ddots \\
& & s_k
\end{pmatrix}
\subset GL(k) 
$$
with the Chern roots of the tautological bundle over $X$, and with 
the Chern classes of the tautological line bundles over $X_S$. At 
fixed points of $X_S^\bA$ 
$$
s_i = a_{\mu(i)}^{-1} \,, \quad
i=1\,\dots\, k\,. 
$$
We choose upper-triangular matrices as $B$. 
Then the special lift $F'$ of the fixed point from 
Section \ref{s_choiceF} is such that 
$$
\cN^\vee|_{F'} = \sum_{i<j} s_j/s_i = \sum_{i<j} a_{\mu(i)}/
a_{\mu(j)} < 0 
$$
which means that the map \eqref{mapmu} is increasing, that is, 
\begin{equation}
1\le \mu(1) < \mu(2) < \dots < \mu(k) \le n \, .\label{mapmu2}
\end{equation}

\subsubsection{}
The polarization $T^{1/2} X_S$ in \eqref{polpull} differs from the standard 
polarization of $X_S$ by 
$$
\delta T^{1/2} = \hbar^{-1} \cN^\vee - \cN \,, 
$$
which by \eqref{dT12} shifts the K\"ahler parameters by 
$2\hbar\rho$, where 
$$
2\rho = (n-1,n-3, \dots, 1-n)
$$
is the sum of positive roots. 

\subsubsection{}
Denote by
$$
\bff_m(s,z)=\prod_{i<m} \vth(sa_i) \,
\frac{\vth(sza_i\hbar^{m-n})}{\vth(z\hbar^{m-n})} \,
\prod_{i>k} \vth(sa_i\hbar)
$$
the function from \eqref{StabPn}. Let $F_\mu \in X^\bA$ be the 
fixed point below the point \eqref{mapmu2}. 
We have
\begin{equation}
\Stab(F_\mu) = \textup{Symm} \,
\frac{\prod_{i=1}^k \bff_{\mu(i)} (s_i, z\hbar^{2\rho_i})}
{\prod_{i<j} \vth(s_i/s_j) \, \vth(s_j/s_i/\hbar)}\label{StabGrass} \,.
\end{equation}
The numerator is the stable envelope on $X_S$, with the shift of 
the K\"ahler variables explained above. One of the terms 
in the denominator is $\Theta(\hbar^{-1}\cN^\vee)$ that comes
from inverting $\bj_{-,*}$.  The other term in the denominator, together with the symmetrization
in the variables $s_i$, 
is the push-forward $\pi_*$. We wrote $\Theta(\cN^\vee)$ with a 
change of sign, so that the sign agrees with polarization. 

Formula \eqref{StabGrass} is the elliptic analog of Proposition
5.3.1 from \cite{Shen}.

\subsubsection{}

When the fixed loci $X_S^\bA$ are not isolated, abelianization
formulas become more involved but can still be made reasonably
explicit, see \cite{SmirHilb}.

\subsection{K-theory limit}

\subsubsection{}

As $q\to 0$, the elliptic curve $E$ converges to the nodal 
elliptic curve with smooth locus isomorphic to $\C^\times$ and 
the group law given by $(x,y)\mapsto x+y+xy$. Correspondingly, 
the scheme $\Ell_\bT(X)$  converges to the spectrum of 
$K_\bT(X) \otimes \C$. In analytic terms, one can see these 
flat limits of schemes
very concretely, as $\Ell_\bT(X)$ are quotients of 
$\textup{Spec} \, K_\bT(X) \otimes \C$ with a fundamental 
domain that grows to become the whole space.

Sections of line bundles on $\Ell_\bT(X)$, which are functions with a
certain automorphy under multiplicative translations by 
$q$ and $e^{2\pi i}$, 
converge to sections of line bundles on $K_\bT(X) \otimes \C$, i.e. to
function that may pick up sign under a  multiplicative translation by 
$e^{2\pi i}$. 
Restricted to $X^\bA$, these become Laurent polynomials on $\bA$ with 
integral or half-integral exponents.

In the $q\to 0$ limit, the effect of shifts by $q$ is replaces by the
growth condition on the these Laurent polynomials at the infinity of
$\bA$. This growth condition, which is the subject of this subsection,
is formulated in terms of the Newton polygons, matching the definition
of the K-theoretic stable envelope, see \cite{Opcmi}.

\subsubsection{}
Recall that the function \eqref{theta_function} satisfies
$$
\vth(e^{2\pi i} x) = - \vth(x) \,.
$$
To avoid half-integer exponents, we define 
$$
\Stab_\bigtriangledown = \left(\det T^{1/2}\right)^{-1/2} \circ
\Stab\circ 
\left(\det T_{X^\bA}^{1/2}\right)^{1/2} \,. 
$$
{}From \eqref{stab_diag} 
we get 
\begin{equation}
\Stab_\bigtriangledown \to (-1)^{\rk \ind} 
\left( \frac{\det N_{<0}}{\det N^{1/2}} \right)^{1/2}  j_*
\pi^*\,, \quad q\to 0\,,
\label{stab_diag_K}
\end{equation}
near the diagonal, 
where $j_*$ and  $\pi^*$ are the corresponding operations in
K-theory. The restriction of $ j_* \pi^*$ back to the fixed locus
is given by tensor product with 
$$
\Ld \, N_{<0}^\vee = \sum_k (-1)^k \Lambda^k N_{<0}^\vee \,.
$$
In particular, the $\bA$-weights in \eqref{stab_diag_K} are the 
same as those occurring in the expansion of $\prod_{w\in N^{1/2}}
(1-w^{-1})$.

\subsubsection{} 

For a given component $F$ of $X^\bA$, we denote 
$$
\Delta_F = \textup{Convex hull} \left( \supp \prod_{w\in N^{1/2}} (1-w^{-1}) 
\right) \subset \bA^\wedge \otimes_\Z \R \,,
$$
where support means the set of weights that occur in the expansion of
the product, that is, the set of exponents that appear in a
multivariate Laurent polynomial.

A convex hull like this is called a Newton polytope and it controls
the growth of the polynomial at the infinity of a torus $\bA$. 
The $q\to 0$ of
elliptic stable envelopes will be characterized by a bound on their
Newton polytopes in equivariant variables, in addition to the same
triangularity and normalization as for the elliptic stable
envelopes. This will precisely recover the definition of the
K-theoreric stable envelopes, see \cite{Opcmi}.

\subsubsection{}

Let $F_1$ and $F_2$ be two components of $X^\bA$ and consider the 
restriction 
\begin{equation}
\Stab_\bigtriangledown \Big|_{F_2 \times F_1}\label{Stab_tr_restriction}
\end{equation}
of the stable envelope of $F_1$ to $F_2$. Suppose that $\ln(q)$ and 
$\ln(z)$ both go infinity so that 
\begin{equation}
- \Re \, \frac{\ln z}{\ln q} \to \cL \in \Pic(X) \otimes_\Z \R
\,, \label{lnzlnq}
\end{equation}
where $\cL$ is generic.

\begin{Proposition}
  If $X$ is a hypertoric variety or a Nakajima quiver variety then 
\begin{equation}
  \supp \lim_{q\to 0} \Stab_\bigtriangledown \Big|_{F_2 \times F_1}
\subset \Delta_{F_2} +  \textup{weight of $\cL\big|_{F_2}$} - \textup{weight of
  $\cL\big|_{F_1}$} \,.  \label{window} 
\end{equation}
\end{Proposition}

\noindent
While the statement should not be specific to quiver varieties, we
prove it here in the same generality as the existence
of elliptic stable envelopes.
Also note the multivariate statement \eqref{window} is equivalent to its
  pull-back by an arbitrary cocharacter $\C^\times \to
  \bA$. Therefore,
  it is enough to prove it when $\rk \bA = 1$.

  \begin{proof}

    For hypertoric varieties, the statement follows at once from
    \eqref{s_theta_prod} and
    \eqref{limit_theta}. For Nakajima quiver varieties, we follow the
    logic of the proof of Theorem \ref{t1}.

    Consider the bottom arrow in the diagram \eqref{diagXs}. In the
    composition of the three maps there, the middle 
    map $\bj_+^*$ is the functorial pull-back in elliptic cohomology
    and it does not change restriction to fixed points. The other two
    maps, give theta-functions in the denominators, as exemplified by
    the formula \eqref{StabGrass}.

    We have
    $$
    x^{-1/2} \vth(x) \to 1 - x^{-1} \,, \quad q\to 0 \,,
    $$
    and thus both maps reduce the Newton polytope by precisely the
    difference between the polarizations of $X$ and $X_S$.  This
    concludes the proof. 
\end{proof}

 We have shown the following 

\begin{Proposition}
In the limit when $q\to 0$ and $\ln z\to \infty$ so that the limit 
\eqref{lnzlnq} is generic, the elliptic stable envelope
$\Stab_\bigtriangledown$ converges to the K-theoretic stable envelope 
with slope $\cL$. 
\end{Proposition}

\section{$R$-matrices} \label{s_R_matr} 

\subsection{Definition} 

\subsubsection{}

Let $\fC_1$ and $\fC_2$ be two chambers in $\Lie \bA$. 
We define 
$$
R_{\fC_2 \leftarrow \fC_1} = \Stab_{\fC_2}^{-1} \circ 
\Stab_{\fC_1} \,,
$$
where some fixed polarization is understood. Instead of the inverse, 
we can use the dual, as in Proposition \ref{p_dual}.
By construction, 
$$
\Theta\left(T^{1/2} X^\bA \right) \otimes 
\tau(\hbar \det \delta \ind)^* \cU
\xrightarrow{\, R_{\fC_2  \leftarrow \fC_1} \,}  \Theta\left(T^{1/2} X^\bA
\right) \otimes \cU \,,
$$
restricted to K\"ahler variables of $X$. Here 
$$
\delta \ind = \ind_2 - \ind_1
$$
where $\ind_1$ and $\ind_2$ 
are the attracting parts of the polarization according to the 
$\fC_1$ and $\fC_2$, respectively. 

\subsubsection{}

Suppose $\fC_1$ and $\fC_2$ are separated by a wall, that is, 
share a codimension one subcone $\fC'$. Using the notations of 
Section \ref{sfC'}, we write $X'=X^{\bA'}$ and note that 
 the decomposition 
$$
N_{X/X^\bA} = N_{X/X'} \oplus N_{X'/X^\bA} 
$$
of normal bundles separates the weights into those that have the same
(respectively, opposite) 
 sign with respect to $\fC_1$ and $\fC_2$. Let 
$$
R_\textup{wall} = R_{\fC_2/\fC'  \leftarrow \fC_1/\fC'} 
$$
be the $R$-matrix for $X'$. Proposition \ref{p_tria} 
implies
\begin{equation}
R_{\fC_2  \leftarrow \fC_1}(z) = R_\textup{wall} (z-\hbar \, \det
\ind_{X',\fC'}) \,,\label{RRwall}
\end{equation}
where 
$$
\ind_{X',\fC'} = \left(T^{1/2} X \big|_{X'} \right)_{>_{\fC'} 0}  \,. 
$$

\subsubsection{}
If $\fC_{n+1} = \fC_1$ then evidently 
\begin{equation}
R_{\fC_{n+1} \leftarrow \fC_n} \, \cdots \, 
\cdots R_{\fC_3  \leftarrow\fC_2}  R_{\fC_2  \leftarrow \fC_1} = 1
\,,\label{braid}
\end{equation}
which written in terms of \eqref{RRwall} is a form of Coxeter
relation for the wall $R$-matrices. 

For example, for $n=2$ we get the unitarity relation 
$$
R_{21}(a^{-1}) \, R(a) = 1
$$
satisfied by wall $R$-matrices, where $a\in E =\cE_\bA$ 
is the equivariant parameter. An example of a braid 
relation is the dynamical Yang-Baxter equation for 
tensor products of Nakajima varieties.

\subsection{Tensor products of Nakajima varieties}\label{s_tensor_prod}

\subsubsection{}

The setup is the same as in 4.1.5 of \cite{MO1}. 
A very special, but very important special case of the 
above construction is 
$$
\bA \subset G_W\,, \quad G_W = \prod GL(W_i)\,,
$$
where $W_i$ are the framing spaces of a Nakajima variety $X$. 
The fixed points $X^\bA$ in this case are products of 
smaller Nakajima varieties associated to the same quiver. 

\subsubsection{}
Concretely, let $a\in \bA \cong \Ct$ act on the framing spaces
with the character 
$$
W_i = a \, W'_i + W''_i \,. 
$$
Over the fixed loci the tautological bundles $\cV_i$ split
$$
\cV_i =  a \cV'_i + \cV''_i 
$$
and the components of the fixed loci are parametrized by the 
corresponding splitting $v=v'+v''$ of the dimension vector. 

\subsubsection{}\label{sindVV} 

Choosing a polarization of the Nakajima variety requires a choice 
of the orientation of the quiver. Let $\bQ$ be the 
adjacency matrix of the oriented quiver. We take 
$$
T^{1/2}X = \sum_i \Hom(W_i,\cV_i) + \sum_{i,j} 
(\bQ_{ij}-\delta_{ij}) \, \Hom(\cV_i,\cV_j) \,.
$$
In principle, since $\bT$ acts on the multiplicity spaces, the 
matrix $\bQ_{ij}$ takes values in $K_\bT(\pt)$. By our 
assumption, $\bA$ acts only on the framing spaces and hence 
$$
\ind X^\bA = \sum_i \Hom(W''_i,\cV'_i) +
\sum_{i,j} 
(\bQ_{ij}-\delta_{ij}) \, \Hom(\cV''_i,\cV'_j)  \,,
$$
assuming that $a$ is an attracting weight.

\subsubsection{}
In particular, 
\begin{equation}
  \label{c1ind}
  c_1(\det \ind) = \lan w'' + (\bQ^t-1) \cdot v'', c_1(\cV') \ran 
+ \lan  (1-\bQ) \cdot v', c_1(\cV'') \ran\,. 
\end{equation}
where 
$$
c_1(\cV) = ( c_1(\cV_1), c_1(\cV_2), \dots)
$$
is a vector of generators of $\Pic(X)$ and
$\lan\,\cdot\,,\,\cdot\,\ran$ is the coordinate pairing. 
We identify $\Pic(X)$ with the coordinate lattice accordingly. 

\subsubsection{}

Now suppose $\bA \cong (\Ct)^3$ acts on $W$ so that 
\begin{equation}
W_i = a_1 W^{(1)}_i +  a_2 W^{(2)}_i +  a_3 W^{(3)}_i\label{aWaW}
\end{equation}
and that the weights $a_1/a_2$ and $a_2/a_3$ are attracting for
$\fC$. In this case, there is a total of $6$ chambers, and going 
around them as in \eqref{braid} using \eqref{RRwall} and 
\eqref{c1ind} we get the dynamical Yang-Baxter equation.

To write it more concisely, introduce the Cartan matrix 
$$
\bC = 2 - \bQ - \bQ^t \,, 
$$
and write $\mu = w - \bC \cdot v$. With dynamical 
variables written additively, the dynamical 
Yang-Baxter equation for the matrix 
$$
\bR(z) = R_\textup{wall}(z+\hbar{(1-\bQ) v}) 
$$
reads as follows. 

\begin{Proposition}
The matrix $\bR(z)$ satisfies the equation (dynamical 
parameters written additively): 
\begin{multline}
  \label{DYB}
\bR_{12}\big(z\big) \, 
  \bR_{13}\big(z-\hbar\, \mu^{(2)} \big)  \, 
\bR_{23}\big(z\big) = \\
=
\bR_{23} \big(z- \hbar\, \mu^{(1)}\big) \, 
 \bR_{13} \big(z\big) \, 
\bR_{12} \big(z- \hbar\, \mu^{(3)}\big) \,. 
\end{multline}
\end{Proposition}

\begin{proof}
For $\bA$ as in \eqref{aWaW}, there are two ways to cross 
from the chamber 
$$
\fC_{123}=\{a_1 > a_2 > a_3\}
$$
 to the opposite chamber  $\fC_{321}$.  
Consider the one that goes through $\fC_{132}$ and $\fC_{312}$, that
is, crosses the walls $a_2=a_3$, $a_1=a_3$, 
and $a_1=a_2$ in this order. 

On the wall $a_2=a_3$, we have, in the notation of Section 
\ref{sindVV}, 
$$
W'=W_1, \quad W''=W_2+W_3\,,
$$
and similarly for $\cV'$ and $\cV''$ and the K\"ahler variables $z$ of 
the matrix \eqref{RRwall} correspond to the bundles $\cV''$. Therefore, 
from \eqref{c1ind} we have 
$$
R_{\fC_{132}  \leftarrow \fC_{123}}(z)  = 
R_\textup{wall}(z-\hbar (1-Q) \bv_1) \,. 
$$
On the next wall $a_1=a_3$, we have 
$$
W'=W_1+W_3,\quad W''=W_2\,, \quad \textup{etc.},
$$
while the K\"ahler variables correspond to $\cV'$. Therefore 
\begin{align*}
  R_{\fC_{312}  \leftarrow \fC_{132}}(z)  &= 
R_\textup{wall}(z-\hbar (w_2+(Q^t-1) \bv_2)) \\
&= 
R_\textup{wall}(z-\hbar \mu^{(2)}-\hbar (1-Q) \bv_2) \,. 
\end{align*}
Similarly, 
$$
  R_{\fC_{321}  \leftarrow \fC_{312}}(z) =
 R_\textup{wall}(z-\hbar (1-Q) \bv_3) \,
$$
and so 
\begin{multline*}
R_{\fC_{321}  \leftarrow \fC_{312}} \, 
R_{\fC_{312}  \leftarrow \fC_{132}} \,  
R_{\fC_{132}  \leftarrow \fC_{123}} \Big|_{z\mapsto
  z+\hbar(1-Q)(\bv_1+\bv_2+\bv_3)}  = \\
\bR_{12}\big(z\big) \, 
  \bR_{13}\big(z-\hbar\, \mu^{(2)} \big)  \, 
\bR_{23}\big(z\big) \,. 
\end{multline*}
Doing the same computation for the other way to cross gives
the equation stated. 
\end{proof}

\subsubsection{Example: $\mathfrak{sl}_2$}

For the $A_1$ quiver we have
$$
\bQ = (0) \,, \quad \bC = (2) \,.
$$
The tensor product $\C^2(a_1) \otimes \C^2(a_2)$ is 
identified with the cohomology of the $w=2$ Nakajima 
varieties 
$$
X = \pt \,\,\, \sqcup \,\,\, T^*\bP^1 \,\,\, \sqcup \,\,\, \pt   \,. 
$$
The nontrivial $2\times 2$ block of this matrix is readily 
computed using the results of Section \ref{seTP}. 

For comparison with that section, one should substitute 
the representation \eqref{eqMW} by the $GL(W) \times GL(1)$ 
module 
$$
\Hom(W,V) = W^\vee \otimes \textup{defining $\C$} \,,
$$
which in practical terms means inverting the equivariant 
parameters for $\bA \subset GL(W)$. With this, one 
computes the nontrivial $2\times 2$ block of the $R$-matrix 
as follows: 
\begin{align}
  \bR &= 
        \begin{pmatrix}
          \vth(u) & \dfrac{\vth(\hbar) \, \vth(z u)}{\vth(z)}
          \vspace{5pt} 
\\
& \vth(\hbar/u)
        \end{pmatrix}^{-1} 
 \begin{pmatrix}
          \vth(u\hbar)  \vspace{5pt}  \\
 \dfrac{\vth(\hbar) \, \vth(z/u)}{\vth(z)}
& \vth(1/u)
        \end{pmatrix} \notag \\
& = 
\frac1{\vth(u/\hbar)} 
\begin{pmatrix}
\dfrac{\vth(z\hbar) \vth(z/\hbar) \vth(u)}{\vth(z)^2} & 
-\dfrac{\vth(\hbar) \vth(z u)}{\vth(z)} \vspace{5pt}  \\
-\dfrac{\vth(\hbar) \vth(z/u)}{\vth(z)}  &  
\vth(u)
\end{pmatrix} \label{Rsl2} 
\end{align}
where $u = a_1/a_2\in \chr(\bA)$ and the top left matrix 
element is simplified using 
$$
\vth(A+B)\, \vth(A-B) \, \vth(C)^2 + \textup{cyclic} = 0 
$$
with $(A,B,C)=(z,u,\hbar)$.

This differs by a gauge transformation from Felder's $R$-matrix
\cite{FelICM}
which reads, in current notation, 
$$
\bR_\textup{standard} = 
\frac{1}
{\vth(\hbar/u)\vth(z)}
\begin{pmatrix}
\vth(z \hbar)\vth(1/u)&  \vth(z u )\vth(\hbar) \\ 
\vth(z/u)\vth(\hbar)&  
\vth(\hbar/z)\vth(u) 
\end{pmatrix}
\,.
$$
%
Recall that gauge transformations of the form 
$$
\begin{pmatrix}
  * & * \\
* & * 
\end{pmatrix}
\mapsto 
\begin{pmatrix}
 f(z) * & * \\
* & f(z)^{-1} * 
\end{pmatrix}\,, 
$$
where $f(z)$ is arbitrary, preserve solutions of \eqref{DYB}. 

\subsubsection{}
Similarly, the tensor square of the defining representation 
of $\mathfrak{gl}_n$ is geometrically realized in a union of 
points and $T^*\bP^1$, thus recovering the basic elliptic 
solution of the dynamical Yang-Baxter equation.

\section{Difference equations}
\label{sdiff}

\subsection{Vertex functions}

\subsubsection{}

Elliptic stable envelopes may be used to determine the monodromy 
of certain difference equations that play a key role in $K$-theoretic 
counting of rational curves in a Nakajima variety, see \cite{Opcmi}. 

In the precise technical sense, this counting refers to computation 
in $K$-theory of the moduli spaces of \emph{quasimaps} from 
$\bP^1$ to $X$. This is closely related to computations 
in 3-dimensional supersymmetric gauge theories, with ${\cal N}=4$ supersymmetry, on
real threefolds of the form 
$$
 \textup{Disk} \, \times \, S^1 \,. 
$$
Nakajima varieties appear in this context as Higgs branches of 
the moduli spaces of supersymmetric vacua. 

In spirit, such curve counting is not far from the subject of quantum K-theory, which is
defined by K-theoretic computation on the moduli spaces
$\Mbar_{0,n}(X)$ and $\Mbar_{0,n}(X\times \bP^1)$ of stable 
pointed rational maps to respective targets. 
Givental and his collaborators studied difference equations 
in quantum K-theory \cites{GL,GT}, and their theory is very 
general --- $X$ can be an arbitrary nonsingular 
algebraic variety, or even an orbifold, as long as rational curves
in $X$ satisfy certain properness assumptions.\footnote{Difference equations arising from supersymmetric gauge theories are actively studied by physicists, see for example \cite{pd}.}

By contrast, the development of the theory explained in \cite{Opcmi} 
crucially uses certain specific features of Nakajima varieties and of 
moduli spaces of quasimaps in order to get a better control over the 
difference equations. Eventually, this allows for an explicit 
identification of difference equations in terms of a 
geometric action of a certain quantum group on $K(X)$, see 
\cite{OS}. 

\subsubsection{}

In this paper, we unwrap the complexities of K-theoretic curve
counting only to the extent required to state and prove our results. 
We denote by 
$$
\QM(X) =\{\textup{stable } f: \bP^1 \dasharrow X \}  \big/
\cong 
$$
the moduli space of stable quasimaps to $X$ 
as defined in \cite{CFKM}.  This is a countable union of 
algebraic varieties indexed by
$$
\deg f \in H_2(X,\Z)_{\textup{effective}} \,.
$$
The standard action of $\Ct$ on $\bP^1$ with 
$$
\left(\bP^1\right)^{\Ct} = \{0,\infty\} 
$$
 induces an action of $\Ct$ on 
$\QM(X)$. We denote by $q$ the weight of $T_0 \bP^1$ and denote 
this $\Ct$ by $\Ct_q$ in what follows. 

\subsubsection{}

We denote by 
$$
\QM(X)_{\textup{nonsing at $\infty$}} \subset \QM(X) 
$$
the open set formed by quasimaps nonsingular at $\infty\in \bP^1$. 
Evaluation at this nonsingular point gives a map 
$$
\ev:  \QM(X)_{\textup{nonsing at $\infty$}} \to X \,. 
$$
This map is not proper, but its restriction to the $\Ct_q$-fixed 
locus is proper. 

\subsubsection{}

The vertex function is a multivariate formal power series 
$$
\Vx   \in (K(X)\otimes \C[\Ct_q])_\text{localized} [[z]] 
$$
defined by 
$$
\Vx = \ev_* \Big(
\textup{symmetrized virtual structure sheaf} \,\cdot \, \label{defXi}
z^{\deg} \Big)\,,
$$
where the pushforward is defined using $\Ct_q$-localization and 
$z^{\deg}$ is an element of the semigroup algebra of
$H_2(X,\Z)_\textup{effective}$. 

\subsubsection{}

To talk about difference equations, we need to introduce two 
transcendental functions. The first function is related to the  
operators 
$$
\ln \big( \cL \otimes \textup{---} \big) \in \End K_\bT(X)_\textup{localized} \otimes 
\left(\Lie \bT\right)^* \,.
$$
Localization and extension of scalars is needed because:
\begin{itemize}
\item[---] generalized eigenspaces of $\cL \otimes$ are $K(F_i)$, where 
$F_i$ are the components of the fixed point locus $X^\bT$, 
\item[---] the eigenvalues are the weights of
  $\cL\big|_{F_i}$, 
\item[---] the logarithms of the eigenvalues lie in $\left(\Lie
    \bT\right)^*$. 
\end{itemize}
We define the map 
$$
\bla: H^2(X,\C) \otimes \Lie \bT  \to \End K(X^\bT) \otimes \C
$$
which extends the map 
$$
\cL  \mapsto \ln (\cL \otimes \textup{---} )   
$$
by linearity to $H^2(X,\C) = \Pic(X) \otimes_\Z \C$. 
This is a close relative of the maps from Section \ref{s_tilde_c} 
and the function 
$$
\be(z) = \exp\left(\frac{\bla(\ln z ,\ln t)}{\ln q}\right)
$$
is a section of $\cU^{-1}$.

\subsubsection{}

The other transcendental function we need is related to 
reciprocal
$$
\phi(x) = \prod_{i\ge 0} (1-q^i x) 
$$
 of the $q$-Gamma function, also known under various other 
names such ``quantum dilogarithm''. It is a half of the theta-function in the 
sense that 
\begin{equation}
  \label{theta_phi}
  \vth(x) = (x^{1/2} -x^{-1/2}) \, \phi(qx) \, \phi(q/x) \,. 
\end{equation}
Given $\cG\in K(X)$ we 
define 
$$
\Phi(\cG) = \prod \phi(\gamma_i) \in K(X)[[q]] 
$$
where $\gamma_i$ are the Chern roots of $\cG$.

\subsubsection{}

We have the following 

\begin{Theorem}[\cites{Opcmi,OS}] \label{t_diff_eq} The $K(X^\bT)$-valued function 
  \begin{equation}
    \label{tVx}
    \tVx = \be(z_\#) \, \Phi((q-\hbar)
    T^{1/2}) \, \Vx  
  \end{equation}
where 
\begin{equation}
  \label{z_sharp}
  z_\# = z (- \hbar^{1/2})^{-\det T^{1/2}}
\end{equation}
generates a holonomic module of rank $\rk K(X)$ under the action of 
$q$-difference operators in both K\"ahler and equivariant variables
$a\in\bA$. This 
$q$-difference module has 
regular singularities in $z$ and $a$ separately. 
\end{Theorem}

The part about regular singularities in $\bA$ is shown in \cite{Opcmi}, 
regular singularities in $z$ follow from \cite{OS}. The existence of 
\emph{some} difference equation in $z$ should be also a 
corollary of Givental's theory. The main point of \cite{OS} is a 
specific representation-theoretic identification of a difference 
connection in the K\"ahler variables. Since \eqref{tVx} has a
prefactor of the very shape discussed in 
Section \ref{s_sep_reg}, the equation 
for $\tVx$ will not be regular jointly in $z$ and $a$. 

If $\dim X^\bT>0$, the subbundle of invariants
$$
I = \left(\hbar T^{1/2}  \big|_{X^\bT} \right)^\bT 
$$
may be nontrivial. The corresponding term 
\begin{equation}
\Phi(-I) = \Phi(-q I) \prod_{\textup{Chern roots $\gamma_i$}}
\frac1{1-\gamma_i} \label{Phi-I}
\end{equation}
in the prefactor of $\eqref{tVx}$ is independent of 
of any equivariant or 
K\"ahler variables and thus, from the point of difference
equations, can be simply dropped. It is convenient, however, 
to keep it and interpret the singular part in \eqref{Phi-I} as a section of an analog of Thom 
bundle $\Theta(-I)$ in K-theory.

\subsection{Example: quasimaps to $T^*\bP(W)$}\label{seQMTP}

We compute the vertex functions in the  simplest example, using
 notations of Section \ref{seTP}.  

\subsubsection{}

In brief, quasimaps to a quotient are sections of 
a bundle of prequotients, up to isomorphism, that is,  
$$
\QM\left(\bP^1 \dasharrow Y/G\right) = 
\Big\{\bP^1 \to Y \times_G \cP \Big\}_\textup{stable} \Big/ \cong
$$
where $\cP$ is a principal $G$-bundle over $\bP^1$. In particular, a principal $GL(1)$-bundle over $\bP^1$ is the 
same as a line bundle $\cL\cong \cO(d)$ where $d=\deg \cL$. 
We define 
$$
\tM = W \otimes \cL \,,
$$
and study sections 
$$
(v,\xi) \in H^0(\bP^1,\tM\oplus\hbar^{-1}\,\tM^*)
$$
satisfying $\lan \xi,v \ran =0$, modulo 
$\Aut(\cL)=\Ct$. The twist by $\hbar$ here is to 
record the equivariance with respect to the group
$$
GL(W) \times \Ct_\hbar \times \Ct_q
$$
where $\Ct_q$ acts by automorphisms of the domain $\bP^1$. 

\subsubsection{}

If $v$ is a nowhere vanishing section, then $(v,\xi)$ defines a map 
$
f: \bP^1 \to X
$
of degree $d=\deg \cL$ such that 
$$
\tM = f^* \cM \,,
$$
where $\cM$ is the bundle on $X$ from Section \ref{s_cM}. 
Quasimaps relax this condition on $v$ 
by allowing $v$ to be nonzero as a section. This is the stability 
condition. The points $p\in\bP^1$ where $v(p)=0$ 
are called the singularities of the quasimap. At all 
nonsingular points, a quasimap has a well-defined value in 
$X$. 

\subsubsection{}

The moduli space of quasimaps has a perfect obstruction theory 
with 
\begin{align}
  T_\vir &= \textup{Deformations} - \textup{Obstructions}\notag\\ 
 &= \Hd(\tM + \hbar^{-1} \tM^* - 
(1+\hbar^{-1}) \End
   \cL) \,,
\label{TvirQM} 
\end{align}
where the term with $\End\cL$ accounts for deformation theory 
 of $\cL$ and the moment 
map equation $\lan \xi,v \ran =0$. Of course, for a line bundle, 
$\End\cL$ is canonically trivial, but in the case of $GL(k)$-quotients
this term will be replaced by the endomorphisms of a rank
$k$ vector bundle --- a nontrivial bundle of rank $k^2$.

\subsubsection{}

Quasimaps fixed under $\bA\times\Ct_q$, where 
$$
\bA=\{\diag(a_1,a_2,\dots)\} \subset GL(W)\,, 
$$ 
are 
the following. Let $F_k = \C e_k \subset \bP(W)$ be a fixed point 
in $X$. We take 
$$
\cL = a_k^{-1} \otimes \cO_{\bP^1}(d[0]) 
$$
that is, the sheaf of functions with pole of order $\le d$ at
$0\in\bP^1$, twisted by a character $a_k^{-1}$ of $\bA$. 
This twisting is because $\bA$-fixed point on a quotient by 
$\Aut(\cL)$ defines a map $\bA\to \Aut(\cL)$. 
The fixed quasimap is the unique $\bA\times\Ct_q$-equivariant 
inclusion 
$$
f: \cO \hookrightarrow \cO(d[0]) \hookrightarrow W \otimes \cL \,. 
$$
It is singular at $0\in\bP^1$ and evaluates to $F_k$ under
$\ev_\infty$. The virtual tangent space at this point is easily computed from
\eqref{TvirQM}. The result is that the virtual tangent 
space, minus $\ev^* TX$, equals 
\begin{multline}
  T_\vir(f) - T_{F_k} X =\\
= (q+\dots+q^d) \sum_i \frac{a_i}{a_k} -
  \hbar^{-1} (1+q^{-1}+\dots+q^{1-d}) \sum_i \frac{a_k}{a_i}
  \,.\label{eq:2}
  \end{multline}

\subsubsection{}
The symmetrized virtual structure sheaf is defined by 
\begin{align}
\tO_\vir &= \cO_\vir \otimes \left( \cK_\vir \, 
\frac{\det \cT^{1/2} \big|_{\infty} }{ \det \cT^{1/2} \big|_{0}  } 
\right)^{1/2} \notag\\ 
&=  q^{-\frac12 \deg\cT^{1/2}} 
\cO_\vir \otimes  \cK^{1/2}_\vir  \,, \label{def_tO}
\end{align}
see \cite{Opcmi}, 
where $\cT^{1/2}$ is the bundle of polarizations. It will be convenient to choose the polarization \emph{opposite}
to the polarization in \eqref{polPn}. This gives  
$$
\deg \cT^{1/2} = - nd\,, \quad n= \dim W \,,
$$
and 
$$
z_\# = z (- \hbar^{1/2})^{-\det T^{1/2}} = (-1)^n \, \hbar^{n/2} \, z \,.
$$
To save on constant factors, we forget the contribution of $\cK_X =
\hbar^{1-n}$. With this convention, the vertex function is the
hypergeometric function 
\begin{align}
  \label{vertex_hyper}
  \Vx\big|_{F_k} &= \sum_{d\ge 0} z^d  \, 
\left(-\frac{q}{\hbar^{1/2}}\right)^{dn} \, \prod_i \frac{(\hbar a_i/a_k)_d}{(q
                    a_i/a_k)_d}   \notag \\
&= \bF \left[ \left. 
\begin{matrix} \hbar a_1/a_k, &\hbar a_2/a_k, &\dots \\
                      q a_1/a_k, &q a_2/a_k, &\dots
\end{matrix} \,\, \right| \, \left(\tfrac{q}{\hbar}\right)^n z_\# \,
                                               \right] \,.
\end{align} 
The notation inside $\bF$ refers to numerators and 
denominators in a hypergeometric series. Readers not 
familiar with these conventions 
may treat the second line in \eqref{vertex_hyper} as the definition of
$\bF$. 
Further 
\begin{equation}
\tVx \big|_{F_k}= e^{-\dfrac{\ln z_\#  \ln a_k}{\ln q}}\, 
\prod_{i\ne k} \frac{\phi(q \hbar^{-1} a_k/a_i)}{\phi(a_k/a_i)} \, 
\bF \left[ \left. 
\begin{matrix} \hbar a_i/a_k \\
                      q a_i/a_k
\end{matrix} \,\, \right| \, \left(\tfrac{q}{\hbar}\right)^n z_\# \,
\right]\label{tVxPn} \,,
\end{equation}
because the weight of $\cO(1)$ on the $k$th fixed point is
$a_k^{-1}$. 

\subsubsection{}

To see directly that the function \eqref{tVxPn} satisfies a 
difference equation in all variables, we can write it as follows 
\begin{equation}
\tVx \big|_{F_k} = \frac1{2\pi i}  \,
\int \frac{ds}{s} \,  e^{\dfrac{\ln z_\# \, \ln s}{\ln q}} 
\Phi'((q-\hbar) \, \textup{Polarization}) \label{VtxInt} \,,
\end{equation}
where
\begin{equation}
\textup{Polarization} = -\frac{1}{\hbar} + \sum_i \frac{1}{\hbar a_i
  s}\,,\label{polprefac}
\end{equation}
the contour of integration enclosed the poles 
$$
s =  \frac{q^d}{a_k} \,, \quad q=0,1,2,\dots \,, 
$$
and prime in \eqref{VtxInt} means we drop the zero 
factor, that is, we define 
$$
\Phi'(1) = \phi(q) \,. 
$$
The integrand in \eqref{VtxInt} obviously satisfies difference 
equations in $s$ and $a_i$. The integration as in \eqref{VtxInt}
acts by Fourier transforms on difference equations, whence 
the conclusion.

Integral representations of vertex functions, of the form \eqref{VtxInt}, are commonplace 
in the literature on supersymmetric gauge theories.
They may be 
interpreted using  heuristic presentation 
\begin{equation}
\QM(X)_\textup{ns at $\infty$} \approx 
\QM(T^*\rM)_\textup{ns at $\infty$}  \rd 
\textup{gauge transformations} \label{QMgauge}
\end{equation}
of quasimaps to an algebraic symplectic reduction $X = T^*\rM \rd  G$, 
as in \eqref{XrdM}. The formula \eqref{QMgauge} lands
on the solid ground of algebraic geometry once one takes
fixed points of the $\Ct_q$-action and this enough to 
produce integral
representations, see the Appendix in 
\cite{AFO}.

While hypergeometric functions such as those in \eqref{vertex_hyper} are well known to have both the series and integral representations, the underlying physics is interesting. The function $\tVx$ computes the partition function of the 3d gauge theory whose Higgs branch is $X$. The series representation of $\tVx$, in terms of summing up quasi maps to $X$, arises when one computes the partition function as a sum over vortex instanton contributions, on the Higgs branch. The same partition function has an integral representation tied to the Coulomb branch of the 3d gauge theory instead (the Chern roots of vector bundles one ends up integrating over are related to the Coulomb moduli). The later is much simpler in physics terms, since it is a result of a one loop computation.\footnote{See \cite{Tudor} for a fairly accessible review, \cite{HC} for more details. See also \cite{WP} for an early manifestation of the same phenomenon.} 



\subsection{Pole subtraction and monodromy} \label{s_ps} 

\subsubsection{}
As already stressed in Section \ref{s_sep_reg}, the difference 
equations satisfied by the vertex functions \eqref{tVx} 
have regular singularities in K\"ahler variables and, 
separately, in equivariant variables $a\in \bA$.  

The solutions $\tVx$ favor K\"ahler variables -- they are given 
by a convergent power series in $z$ and are therefore 
holomorphic in some punctured neighborhood of the point 
$z=0$. Here and below, we allow functions not to be single-valued
in this punctured neighborhood, as their monodromy around 
the origin is an integral part of the story. 

Instead of solutions holomorphic in $z\to 0$, one can fix a point 
$0=0_\fC \in \overline{\bA}$ as in \eqref{0inbA} and ask for a 
basis $\Vx_\fC$ of solutions holomorphic as $a\to 0$. 
We call the transition matrix $\fP$  between these two 
bases of solutions the \emph{pole subtraction matrix}. 

Evidently, matrix elements of $\fP$ are sections of certain 
line bundles on $\cE_{\Pic(X)} \times \cE_\bA$ as functions of 
$z$ and $a$ invariant under $q^{\Pic(X) \oplus \cochar \bA}$, 
whose transformation under 
$$
2\pi i (\Pic(X) \oplus \cochar \bA) \hookrightarrow 
H^2(X,\C) \oplus \Lie \bA  \owns (\ln z,\ln a) 
$$
is dictated by the exponential in \eqref{tVx}. 

\subsubsection{}

A choice of $\fC$ defines a partial order of the components $F_i$ of $X^\bT$, 
called the \emph{ample} partial order in \cite{MO1}. This order is by the order 
of vanishing of the $\bA$-weight of an ample bundle as $a\to 0$. 
The order is coarser than the order from Section \ref{s_order_F}. 

\begin{Proposition}\label{p_growth_order} 
The matrix $\fP$ is block-triangular with respect to the ordering of 
components of $X^\bT$ which is 
\underline{opposite}
to ample. 
\end{Proposition}

\begin{proof}
With the exception of the exponential prefactor, all terms in
\eqref{tVx} grow at most polynomially as $(z,a)\to (0,0)$. Indeed,
for instance 
\begin{equation}
\frac{\phi(b/a)}{\phi(c/a)} \sim \const \, a^{\log_q
    (b/c)} \label{asymp_phi}
\end{equation}
along any geometric progression of the form $a=q^n x$, $x\ne q^n c$, 
because of the difference equation that this function satisfies. 

The exponential prefactor thus determines the rate of growth 
of the solution as $(z,a)\to (0,0)$ and it filters the space of 
solutions by subspaces where this rate of growth is bounded above. 
This means bounds from above on the $\bA$-weight of an ample 
bundle, whence the conclusion. 
\end{proof}

\subsubsection{}
The matrix $\fP$ may, in principle, be computed algorithmically
as follows. We first transform the basis $\tVx$ into a basis 
of solutions which has no poles in the $\Phi$-prefactor and for
which the pole subtraction matrix may be chosen unitriangular 
with respect to the ordering in Proposition \ref{p_growth_order}.

This first step is accomplished by exchanging the terms 
$$
\frac{\phi(q \hbar^{-1} a^{-1})}{\phi(a^{-1})} 
\,\, \longleftrightarrow \,\,
\frac{\phi(q a)}{\phi (\hbar a)} \, a^{1-\log_q \hbar} 
$$
which solve the same difference equation in $a$. This 
exchange may be interpreted as switching between 
$a$ and the opposite weight $\hbar^{-1} a^{-1}$ in a 
polarization. 

On the next step, one considers the poles in the series in \eqref{tVx} 
at divisors of the form $w=q^n$, 
where $w$ is a weight of $\bT$, positive on $\fC$ and $n\gg 0$. 
One observes that such poles occur only in terms of degree $\ge
d(n)$, where $d(n)$ grows linearly with $n$. The singular terms 
in the $w\to q^n$ expansion solve the same $q$-difference equation 
in $z$, specialized at $w=q^n$. 
However, their order of vanishing as $z\to 0$ is higher by at least
$d(n)$. Therefore, they are linear combinations of slower growing 
solutions and adding a suitable linear combination of 
slower growing solutions makes the given solution pole-free at 
$w=q^n$, $n\gg 0$. The matrix coefficients of $\fP$ are 
thus determined as the unique section of a certain line nontrivial 
line bundle on $\cE_A$ with prescribed singularities.

\subsubsection{} 

The main point of this Section is that elliptic stable envelopes give
the pole subtraction matrix $\fP$. To put elliptic cohomology and 
K-theory on the same footing, we will work in analytic completion 
of K-theory, that is, with analytic functions on 
$$
\fK_T(X) = \Spec
K_\bT(X)\otimes \C \,. 
$$
Localization introduces meromorphic functions with poles in 
specified locations. 

For consistency with elliptic cohomology, we think of pushforward 
under $f: X\to Y$ as 
$$
f_*:  \Thom(-N_f) \to \cO_{\fK(Y)} 
$$
even though in K-theory we have Thom isomorphism
$\Thom(-N_f)\cong\cO_{\fK(X)}$.  Then, in particular, the transpose 
of a map $g$ in K-theory acquires an additional twist by the Thom 
class of the tangent bundle, as in Section \ref{s_dual}. 

The transpose in K theory and elliptic cohomology are related by 
$$
g^{\textup{transpose,K}} = \bph^{-1} \, g^\vee \bph
$$
where $g^\vee$ is the elliptic transpose as in Section \ref{s_dual}
and 
\begin{equation}
\bph = \cK_X^{1/2} \, \Phi(q(TX+T^\vee X)) \,. \label{def_bph}
\end{equation}

\subsubsection{}

Define 
\begin{align}
\Stab^\# &=
           \left(\Stab_{-\fC,T^{1/2}}(z_\#^{-1})\right)^{\textup{transpose,K}}
  \label{Stab_sharp2}\\
 &= \bph^{-1} \Stab_{\fC,T^{1/2}_\textup{opp}}(z_\#)^{-1} \label{Stab_sharp} \bph
\end{align}
where $\fC$ is the cone corresponding to the point $0=0_\fC \in
\overline{\bA}$ at which we subtract the poles, and $z_\#$ was 
defined in \eqref{z_sharp}. The equality in \eqref{Stab_sharp} is 
the content of Proposition \ref{p_dual}. 

We define 
\begin{equation}
z_{\#,\fC} =  z_{\#} \, \hbar^{\det T^{1/2}_{<0}} \,, \label{zsfC}
\end{equation}
where $T^{1/2}_{<0}$ stands for the repelling part of the polarization
with respect to the cone $\fC$. The difference with standard shift by 
$\hbar^{\ind} = \hbar^{\det T^{1/2}_{>0}}$ from the earlier sections
is due to the appearance of the opposite cone in \eqref{Stab_sharp2}, 
or the opposite polarization in \eqref{Stab_sharp}. 

{}From the definition of elliptic stable envelopes the operator 
\begin{equation}
  \label{defbS}
 \fP_{\fC} = \be(z_{\#,\fC}) \,  \Theta(T^{1/2} X^\bA) \, \Stab^\#
\, \Theta(T^{1/2})^{-1}  \,  \be(z_{\#})^{-1} 
\end{equation}
commutes with shifts of both $z$ and $a$. 
Here $T^{1/2} X^\bA$ is the polarization of $X^\bA$ induced by
$T^{1/2}$. The term $\Theta(T^{1/2} X^\bA)$ is independent of 
the variables in $\bA$ and is inserted here only for the event that a 
larger torus is acting preserving the symplectic form.

As we will see, the operator \eqref{defbS} 
subtracts the poles as follows 

\begin{Theorem}\label{t_ps} 
The function  
  \begin{equation}
    \label{ThVx}
    \Vx_\fC = \fP_\fC \,  
  \tVx  
  \end{equation}
solves the same scalar difference equations as \eqref{tVx} and is 
holomorphic in a punctured neighborhood of $0_\fC \in
\overline{\bA}$. Equivalently, the function 
  \begin{equation}
    \label{t_ps_st}
    \Stab^\# 
    \, \frac{(\det T^{1/2})^{-1/2}}{\Phi(T^\vee)} \,  
  \Vx  
  \end{equation}
has no poles in $a$ as $a\to 0_\fC$. 
The exact same pole cancellation property is 
true for vertex functions with \underline{descendents}. 
\end{Theorem}

The proof of Theorem \ref{t_ps} will be given in Section 
\ref{s_t_ps}. 
For a discussion of vertex functions with descendents, see \cite{Opcmi}. 

The equivalence of pole cancellation in \eqref{ThVx} and
\eqref{t_ps_st} follows from 
\begin{equation}
\frac{\Phi((q-\hbar) T^{1/2})}{\Theta(T^{1/2})} = 
\frac{(\det T^{1/2})^{-1/2}}{\Phi(T^\vee)}  \,. 
\label{phiThe_phi}
\end{equation}
As before, the 
$\Phi(T^\vee)$-term here is interpreted using
$$
\frac{\Phi(q T^\vee)}{\Phi(T^\vee)} : \cO_\fK \xrightarrow{\,\sim\,}
\Thom(-T)  \,. 
$$

\subsubsection{}

To check the truth of Theorem \ref{t_ps} in the simplest example of 
$T^*\bP^{n-1}$, we change the integrand in \eqref{VtxInt} as in 
\eqref{t_ps_st}, which gives 
\begin{equation}
\Vx_{\fC} \big|_{F_k} = \frac1{2\pi i}  \, e^{\dfrac{\bla(\ln z_{\#,\fC},\ln
    t)}{\ln q}} 
\int_{|s|\approx 1}\frac{ds}{s} \,  \frac{(\det T^{1/2})^{-1/2}
  \,\Stab_\#}{\Phi'(T^\vee)} 
\label{ThVxInt} \,,
\end{equation}
where
$$
T^{1/2} =  -\frac{1}{\hbar} + \sum_i \frac{1}{\hbar a_i s}\,,
$$
is the polarization on the prequotient as in \eqref{polprefac}, 
$$
\Stab_\# = \prod_{i<k} \vth(s a_i \hbar) \,
\frac{\vth(s z_\#^{-1} a_k \hbar^{n-k+1})}{\vth(z_\#^{-1} \hbar^{n-k+1} )} \,
\prod_{i>k} \vth(s a_i)
$$
is the replacement of \eqref{StabPn} for the opposite polarization and 
the opposite ordering, and the difference between $\Phi(qT^\vee)$ and 
$\Phi(T^\vee)$ corresponds to the pushforward $K_\bT(X) \to
K_\bT(\pt)$.  The contour of integration in \eqref{ThVxInt} is roughly 
the compact torus in $\Ct$, it will be specified precisely in a
moment. 

Neglecting irrelevant factors, the integrand in \eqref{ThVxInt}
simplifies to 
\begin{multline}
  \label{eq:4}
  \textup{Integrand} \propto \prod_{i<k} \frac{\phi(\frac q\hbar 
a_i^{-1} s^{-1})} {\phi(a_i^{-1} s^{-1})} 
\prod_{i>k} s a_i \frac{\phi(q a_i s)} {\phi(\hbar a_i s)} \,  \times \\
(s a_k)^{1/2} \, \frac{\vth(s z_\#^{-1} a_k
  \hbar^{n-k+1})}{\vth(z_\#^{-1} \hbar^{n-k+1} )}
\frac{1}{\phi(a_k^{-1} s^{-1}) \phi(\hbar a_k s)} \,,
\end{multline}
which has poles at 
\begin{alignat}{3} 
  s&=q^n a_i^{-1}\,, \qquad && i\le k \,, \quad
  &&n=0,1,\dots \label{pl1} \\
 s&=q^{-n} \hbar^{-1} a_i^{-1}\,, \qquad && i\ge k \,, \quad
 &&n=0,1,\dots\,. \label{pl2}
\end{alignat}
The poles in \eqref{pl1} accumulate to $s=0$ while the poles in
\eqref{pl2} accumulate to $s=\infty$. The contour of the integration
in \eqref{ThVxInt} separates \eqref{pl1}  from \eqref{pl2}. 

A contour integral becomes singular when the poles on opposite sides 
of the contour coalesce, which means 
$$
\frac{q^{n_1}}{a_i} = \frac{q^{-n_2}}{\hbar \, a_j} 
\quad \Rightarrow \quad \frac{a_i}{a_j} = \hbar \, q^{n_1+n_2} \,,
$$
with $i \le k \le j$ and $n_1,n_2 \ge 0$. By our conventions $a\to
0_\fC$ means that $a_i/a_j \to \infty$ for $i<j$ and so the integral 
is pole-free in this region. 

The general statement of Theorem \ref{t_ps} may, in principle, be
approached 
from a similar angle. We will use other, more geometric, tools in what
follows.

\subsubsection{} 

Theorem \ref{t_ps} determines the monodromy of the difference 
equations in $a$ and constraints the monodromy of the 
difference equations in K\"ahler variables in terms of the analogous 
monodromy for $X^\bA$. The following is immediate 

\begin{Corollary}\label{corMonqKZ}
Let $\Vx_{\fC_1}$ and $\Vx_{\fC_2}$ be functions
  \eqref{ThVx} for two different points 
$$
0_{\fC_1}, 0_{\fC_2} \in \overline{\bA} 
$$
then 
$$
\Vx_{\fC_2} = \fP_{\fC_2} \,  \fP_{\fC_1}^{-1} \, \Vx_{\fC_1} 
$$
where the matrices
\begin{equation}
  \label{R_gauge}
\fP_{\fC_2} \,  \fP_{\fC_1}^{-1} = U_2 \,\, \left(
R_{\fC_2 \leftarrow \fC_1, T^{1/2}_\textup{opp}} (z_\#)  \right) \,\, U_1^{-1} 
\end{equation}
with $U_i = \be(z_{\#,\fC_i})   
\frac{\Theta(T^{1/2} X^\bA)}{\bph(X^\bA)}$ 
are gauge transformations of the elliptic $R$-matrices 
from Section \ref{s_R_matr}. 
\end{Corollary}

\subsubsection{}

Let $\tVx_{X^\bA}$ denote the function \eqref{tVx} for the $\bA$-fixed
locus. Note, in particular, that the prefactor in $\tVx_{X^\bA}$ has
the form $\be(z_{\#,X^\bA})$ with 
\begin{equation}
  \label{zsXA}
  z_{\#,X^\bA} = z_\# (-\hbar^{1/2})^{\det N^{1/2}} \,,
\end{equation}
where $N^{1/2}$ is the normal part of the polarization. 
We have the following 

\begin{Proposition}\label{p_a0} 
As $a\to 0_\fC$, we have
\begin{equation}
\Vx_\fC \sim \dots 
\tVx_{X^\bA} \Big|_{z \mapsto z \left( - \hbar^{-1/2}\right)^{\deg
    N_{>0}} \, q^{-\deg T^{1/2}_{<0}}}  \,,
\label{asy_Va0} 
\end{equation}
where dots stand for a factor of the form $\pm \hbar^{k/2}$ which 
depends on the component of the fixed locus. 
\end{Proposition}

In the proof, we will need some generalities about 
the behavior of the solutions of $q$-difference equations at 
a regular singular point. 

Let $f_1(x)$ be analytic 
in a universal cover of a punctured neighborhood of $x=0$
and solve a regular $q$-difference equation. A scalar equation 
may be replaced by an equivalent first order vector equation
\begin{equation}
f(qx) = M(x) f(x) \,, \quad M(0) \in GL(n,\C) \label{fqMf}
\end{equation}
for a vector-valued function, in which $f_1(x)$ is the first entry. 
In fact, from both geometric and representation-theoretic 
viewpoints, it is the matrix equations that arise naturally \cite{Opcmi}. 
We can make a further nonresonance assumption that $\mu_i/\mu_j \ne 
q^k$, $k\ne 0$, where $\mu_i$ are the eigenvalues of $M(0)$. 
This assumption is satisfied in our case. In fact, for the shifts 
of equivariant variables the operator $M(0)$ is semisimple, with 
eigenspaces given by 
\begin{equation}
K_\bT (X^\bA) = \bigoplus_{\textup{components $F$ of $X^\bA$}} \,
K_\bT(F)\label{KXAF}
\end{equation}
and exponents at $a\to 0_\fC$ are the following 
\begin{equation}
  \label{exponents_a}
  \textup{exponents} = \left\{ \left(z_{\#,\fC} \, q^{-\det T^{1/2}_{<0}}
  \right)^{\bla(\,\cdot\,,\sigma)}\right\} \,,
\end{equation}
where $\sigma\in
\cochar \bA$ is the direction of the shift. This is clear from
\eqref{tVx} and \eqref{asymp_phi}. 

With this nonresonance assumption, the 
solution to \eqref{fqMf} has the form exemplified by 
\eqref{tVx} 
\begin{equation}
  \label{seriessol}
  f(x) = e^{\displaystyle \ln M(0) \, \frac{\ln x}{\ln q}}  \times
  \textup{convergent power series in $x$} \,,
\end{equation}
see for example \cite{Gaz} for a modern discussion of this classical
result. In particular, generalized eigenspaces of $M(0)$, given by 
\eqref{KXAF} in our situation, correspond 
to generalized eigenspaces of the monodromy around the origin.
By construction, the solutions \eqref{ThVx} are precisely grouped according 
to the components of $X^\bA$, which gives the following 

\begin{Lemma}\label{l_a0} 
The function $\be\left(z_{\#,\fC} \, q^{-\det T^{1/2}_{<0}}\right)^{-1} 
\Vx_{\fC}$ is analytic in $a$ near $0_\fC$. 
\end{Lemma}

\begin{proof}[Proof of Proposition \ref{p_a0}]
As explained in the proof of Proposition \ref{p_growth_order}, the 
solutions $\Vx_\fC$ are linear combinations of solutions with different order of 
growth as $(z,a) \to (0,0)$.  The function in Lemma \ref{l_a0} is 
meromorphic in $z$, it is therefore uniquely determined by its 
values for $|z| < \varepsilon$ for any $\varepsilon$. 
We are free to choose $\varepsilon$ so 
small that the diagonal term in stable envelopes dominates the 
$a\to 0$ asymptotics. 

On the diagonal, we weights along the fixed locus cancel out 
in the operator \eqref{defbS}. For the weights normal to the 
fixed locus we get, using \eqref{phiThe_phi}, the following 
contribution 
$$
(\det N^{1/2})^{-1/2}  \prod_{w \in N_{>0}}
\frac{\vth(w)}{\phi(w^{-1}) \phi(\hbar w)} \sim 
\hbar^{\dots} \, 
\left(\det T^{1/2}_{<0}\right)^{-1} \,, \quad  a\to 0_\fC \,,
$$
and this is absorbed by the power of $a$ that comes from 
the $q$-shift in the $\be$-factor in Lemma \ref{l_a0}. 

The asymptotics of the vertex function $\Vx$ as $a\to 0$ is 
given by Section 7.3 of \cite{Opcmi} as follows 
$$
\Vx \to \cdots \Vx_{X^\bA} \Big|_{z \mapsto z \left( - \hbar^{-1/2}\right)^{\deg
    N_{>0}} \, q^{-\deg T^{1/2}_{<0}}}  \,. 
$$
Since 
$$
z_{\#,\fC} = z_{\#} \, h^{\det T^{1/2}_{<0}}  = z_{\#,X^\bA} 
(-\hbar^{-1/2})^{\det N_{>0}}  \,, 
$$
the proposition follows. 
\end{proof}

\subsubsection{}

Proposition \ref{p_a0} constrains the monodromy of the difference 
equations in K\"ahler variables as follows. 

All flops of Nakajima 
varieties are the same Nakajima varieties with a different choice
of stability parameter. We identify $\Pic(X)$ with $\Pic(X_\flop)$ 
by sending $\det V_i$ to the same line bundle on $X_\flop$. 
The  decomposition of $H^2(X,\R)$ into ample cones of 
different flops is the chamber decomposition for a certain finite 
collection of rational hyperplanes. The toric variety associated
to the fan of ample cone is the K\"ahler moduli space, the base of the difference 
equations in the K\"ahler variables. Its fixed points, denoted by $0_X$,  $0_{X_\flop},$..., 
correspond to all possible flops of $X$. 

One can compare the solution $\tVx_X$ which is holomorphic 
as $z\to 0_X$ with the solution $\tVx_{X_\flop}$ holomorphic 
as $z\to 0_{X_\flop}$. The comparison is given by the 
monodromy operator
$$
\tVx_{X_\flop} = \Mon_{X_\flop  \leftarrow  X} \, \tVx_X \,,
$$
which is meromorphic and invariant under the shifts of 
variables by $q$. 

\begin{Proposition}
For any chamber $\fC$, the diagram 
  \begin{equation}
\xymatrix{
\Vx_{X,\fC}  
\ar[dd]_{\Mon_{X^\bA_\flop  \leftarrow  X^\bA} 
\Big|_{z\mapsto z \left( - \hbar^{-1/2}\right)^{\deg
    N_{>0}}} }
&& 
\ar[ll]_{\fP_{\fC}}  
\tVx_X \ar[dd]^{\Mon_{X_\flop  \leftarrow  X}}\\ \\
\Vx_{X_\flop,\fC}  && \ar[ll]_{\fP_{\fC,\flop}}  \tVx_{X_\flop}
}
 \label{Mon_fl} 
\end{equation}
commutes.  
\end{Proposition}

\begin{proof}
The connection matrix between the solutions $\Vx_{X,\fC}$ and 
$\Vx_{X_\flop,\fC}$ is invariant under $q$-shifts of $a$ and has 
a limit as $a\to 0_\fC$ by Proposition \ref{p_a0}. Therefore, it 
is constant equal to its value at $a\to 0_\fC$. Since it is
$q$-periodic in $z$, we may ignore shift by $q$ in \eqref{asy_Va0}. 
\end{proof}

\subsection{Proof of Theorem \ref{t_ps}}  \label{s_t_ps} 

\subsubsection{}\label{s_pole_1}

Consider the pole of the vertex function at an irreducible 
divisor of the 
form 
\begin{equation}
w^l q^m = \zeta\,, \quad \zeta^n=1\,, \quad w \in
\left(\bT^\wedge \label{div1} \right)_\textup{indivisible}\,,
\end{equation}
where we can assume that 
$$
\gcd(l,m) =1\,, \quad m>0 \,, \quad \textup{order} \, (\zeta)=n \,.
$$
Consider the subgroup 
$$
\bT  \supset \bT' = \Ker w^{nl} \cong (\Ct)^{\rk \bT-1} \times \mu_{nl}
$$
where $\mu_{nl} \subset \Ct$ is the group of roots of unity. 
Equivariant localization on $\QM(X)$ with respect to $\bT'$, introduces 
poles at divisors of the form 
\begin{equation}
w' q^{m'} =1 \label{div2}
\end{equation}
where $w'$ is a weight of $\bT$ which is nontrivial on $\bT'$. 

Since the poles \eqref{div2} are distinct from \eqref{div1} and 
localization contributes 
$$
\begin{pmatrix}
\textup{ virtual normal } \\
\textup{ bundle terms }
\end{pmatrix} \in 
K_\bT(\QM(X))_\textup{localized at \eqref{div2}} \,, 
$$
we can replace the $\bT$ action on 
$\QM(X)$ by the action of 
$$
\bT_\textup{new} = \bT/\bT'  \cong \C^\times  \owns t 
$$
on 
$$
\QM(X)^{\bT'} = \QM\left(X^{\bT'} \right) = \QM(X_\textup{new}) \,. 
$$
in the analysis that follows. By construction, the 
weight $w^{nl}$ becomes the coordinate on $\bT_\textup{new}$. 
Therefore,  the poles of interest now have the form 
\begin{equation}
t \, q^{m} =1 \label{div3}  \,, \quad m=m_\textup{new} =
n m_\textup{old} > 0 \,. 
\end{equation}

\subsubsection{}
Now suppose $\bT=\Ct_t$ and let 
$$
\bT'' \subset \Ct_t \times \Ct_q
$$
be the subtorus defined by \eqref{div3}. Again, equivariant localization 
with respect to $\bT''$ introduces poles that are distinct from 
\eqref{div3}, therefore it is important to understand the
$\bT''$-fixed loci in the moduli spaces of quasimaps. 

Recall that vertex functions are computed using $\Ct_q$-equivariant
localization. As a scheme,  $\Ct_q$-fixed loci in $\QM(X)_\textup{nonsing at
$\infty$}$ are identical to 
$\Ct_q$-fixed loci in the moduli spaces of twisted quasimaps, 
see \cite{Opcmi}. There is a small but important difference in their
obstruction theory, see below. 

Quasimaps from $\bP^1$ to a GIT-quotient are 
sections of a bundle of prequotients, and one can twist that 
bundle further by using a homomorphism 
$$
\Ct_q \to \Aut(X) 
$$
as a clutching function. This has the effect of allowing $\Ct_q$ to act 
nontrivially in the fibers over $0,\infty\in\bP^1$.  In particular, 
there is a unique twist such that 
\begin{itemize}
\item[---] $\Ct_q$ acts trivially in the fiber over $0$, 
\item[---] $\bT''$ acts trivially in the fiber over $\infty$. 
\end{itemize}
We denote by $\QM(X)_\textup{tw}$ the corresponding moduli space. 

\subsubsection{}
We have the following 

\begin{Theorem}\label{t_proper} 
The map 
\begin{equation}
\ev_\infty: \QM(X)_\textup{tw, nonsing at $\infty$}^{\bT''} \to
X 
\label{evtw}
\end{equation}
is proper for quasimaps of fixed degree and its image lies in 
\begin{equation}
  \label{X_rep}
X_\textup{repelling} = \{x, \lim_{t\to\infty} tx \,\,
\textup{exists}\} \,. 
\end{equation}
\end{Theorem}

\begin{proof}
Singularities of a $\bT''$-fixed quasimap form a $\bT''$-invariant 
finite subset of $\bP^1\setminus\{\infty\}$ and, therefore, they are all confined to the
origin. Thus, away from the origin, we have a $\bT''$-invariant 
parametrized curve in $X$. In local coordinates at $\infty\in\bP^1$ 
this curve is constant, and therefore uniquely determined
by the point in which it meets the fiber at infinity. Given this curve
and degree of the quasimap, 
all possible singularities at the origin form a proper set, whence  
the properness of the map \eqref{evtw}. 

We have proper maps 
$$
X \xrightarrow{\,\,\pi\,\,}  X_0 \hookrightarrow V
$$
where $X_0$ is the affine quotient, that is, the spectrum of the 
algebra of $G$-invariants, and the $\bT$-equivariant embedding 
$X_0\to V$ is obtained from a choice of generators of this algebra. 
This induces proper maps 
$$
\QM(X)_\textup{tw} \to \QM(X_0)_\textup{tw} \hookrightarrow 
H^0(\bP^1,\cV) \,. 
$$
We can split $V$ and $\cV$ by their $t$-weights  
$$
V = V_{\le 0}  \oplus V_{>0}
$$
and then 
$$
X_\textup{repelling} = \pi^{-1}\left( V_{\le 0}\right)
$$
while $\cV_{>0}$ consists of line bundles of negative degree and hence 
has no sections. This shows the evaluation map lands in
$X_\textup{repelling}$. 
\end{proof}

\subsubsection{}

In a discussion of regularity of  functions along a divisor,
it is natural to pass to completion of K-theories at that divisor. 
For example, we may consider 
$$
K_{\Ct_t \times\Ct_q}(\pt)_\textup{completed at \eqref{div3}} 
 = \Q(q)[[t-q^{-m}]] \,. 
$$
This is the completion of rational functions in $t$ and $q$ regular at
$t q^m =1$ in the topology of formal power series in $t-q^{-m}$. 

Further, since elliptic cohomology classes give elements of
$K_\bT(X)[[q]]$, it is natural to additionally pass to the completion 
of the local ring at $q=0$ and define 
\begin{equation} \label{Khat} 
  \widehat{K}_{\Ct_t \times\Ct_q}(\pt)  = \Q((q))[[t-q^{-m}]] \,. 
\end{equation}
It is important to note the order of completions. For the
opposite order, we have 
$$
\frac{1}{1-q^m t} = \sum_{k\ge 0} t^k q^{km} \in \Q[[t-q^{-m}]][[q]]\,,
$$
and so it is meaningless to talk about the order of pole at
$t=q^{-m}$. 

\subsubsection{}

To return to the general setup of Section \ref{s_pole_1}, 
we consider an exact sequence
\begin{equation}
  \label{exTT}
  1 \to \bT' \to \bT \xrightarrow{\,\, t \, \, }  \Ct_t \to 1
\end{equation}
and the corresponding action of $\Ct_t$ on 
$$
X' = X^{\bT'} \,. 
$$
We are interested in poles at the components of 
\begin{equation}
  \label{div4}
  t q^m = 1 \,. 
\end{equation}
Associated to 
this divisor \eqref{div4}, there is the moduli space 
$$
\tQM' = \QM(X')_\textup{tw, nonsingular at $\infty$}
$$
of twisted quasimaps to $X'$, with an 
evaluation map 
\begin{equation}
\ev_\tw^*: \widehat{K}_{\Ct_t \times\Ct_q}(X')_\textup{completed at $t=1$} 
\to \widehat{K}_{\Ct_t \times\Ct_q}\left(\tQM'\right)_
\textup{completed at $q^m t =1$} \label{evK}
\end{equation}
covering the homomorphism of tori 
$$
(t,q)\mapsto (t q^m,q) \,. 
$$
Here the completions of K-theories are as in \eqref{Khat}, with a 
difference of divisors indicated. 

Since $X^\bT=(X')^{\Ct}$ is proper, attracting and repelling sets 
intersect properly. Therefore, from Theorem \ref{t_proper} we 
deduce the following 

\begin{Corollary}
Let 
\begin{equation}
  \cF\in \widehat{K}_{\Ct_t\times\Ct_q} (X')\,,
 \quad
\cG\in \widehat{K}_{\Ct_t\times\Ct_q}(\tQM')
  \end{equation}
be such that 
$$
\supp \cF \subset X'_\textup{attracting} = \{x, \lim_{t\to 0} tx \,\,
\textup{exists}\} \,. 
$$
Then 
\begin{equation}
\chi(\tQM', \cG \otimes \ev_\tw^*
\cF)\in \widehat{K}_{\bT\times\Ct_q}(\pt) \label{chiGF} \,. 
\end{equation}
\end{Corollary}

\subsubsection{} 

Our next goal is to find $\cF$ and $\cG$ so that the
$\Ct_q$-equivariant localization of \eqref{chiGF} reproduces
\eqref{t_ps_st}. We first consider the case when there is a map 
$$
\sigma: \Ct_q \to \bA \subset \bT 
$$
such that 
\begin{equation}
t(\sigma(q)) =  q^m  \,.\label{tsigma}
\end{equation}
Using $\sigma$, we can twist quasimaps to $X$ so that 
$$
\tQM^{\bT'} = \tQM' \,.
$$
And we take in \eqref{chiGF} the restrictions of 
\begin{align}
\cG & = \tO_\vir \otimes \textup{Tautological at $0$} \,, \label{G1} \\
\cF & = \frac{\Stab_\# \otimes \left(\det T^{1/2}\right)^{-1/2}} 
{\Phi(q T^{\vee})} \label{F1} \,,
\end{align}
where the tautological term in \eqref{G1} 
denotes an arbitrary Schur functor of the fibers of 
the tautological bundles at the origin in the domain of the 
quasimap. (For example, this term can be the identity.) 
Note that the twist by the square root of 
$\det T^{1/2}$ makes stable envelope an element of 
$K_\bT(X)[[q]]$ and that the fraction in \eqref{F1} has no pole at
$t=1$. 

\subsubsection{}\label{s_twist_comp1}

We now compare the $\Ct_q$-equivariant localization of \eqref{chiGF} 
with the corresponding computations for untwisted quasimaps. 
The contributions of $\cG$ to two localization formulas are almost 
identical, the only difference comes from 
$$
T_{\vir,\tw} - T_\vir = \frac{(TX)_\tw - TX} {1-q}  \,, 
$$
where $(TX)_\tw$ is the tangent bundle of $X$ with the action of
$\bT\times \Ct_q$ induced by $\sigma$. This gives 
$$
\cO_{\vir,\tw} = \cO_{\vir} \, \Phi\left( q T^\vee_\tw - q T^\vee \right)\,,
$$
which means 
$$
\cO_{\vir,\tw} \otimes \ev_\tw^*\left(\frac1{\Phi(q T^{\vee})}\right) 
\,\, \textup{is twist-invariant} \,. 
$$

\subsubsection{}\label{s_twist_comp2}

The difference between $\cO_{\vir,\tw}$ further includes the weight of 
$\cK_\vir^{1/2}$ and the contribution from the polarization. 
For the latter we have the evident relation 
\begin{equation}
  \label{pol_twist}
  \left(\frac{\det T^{1/2}_{\infty,\tw}}{\det T^{1/2}_0}\right)^{1/2}
  \otimes  \ev_\tw^*\left(\det T^{1/2}\right)^{-1/2} 
\,\, \textup{is twist-invariant} \,. 
\end{equation}
For the former, we have the following 

\begin{Lemma}
  \begin{equation}
    \label{KvirKvir}
    \frac{\cK^{1/2}_{\vir,\tw}}{\cK^{1/2}_{\vir}} = 
\left(-\hbar^{1/2}\right)^{\lan \det T^{1/2}, \sigma\ran} \, 
\frac{\Theta(T^{1/2})}{\Theta(T^{1/2}_\tw)} \,. 
  \end{equation}
\end{Lemma}

\begin{proof}
 This is equivalent to the identity
$$
\det\left(\frac{aq^k + \frac{1}{\hbar a q^k} - a - \frac1{\hbar
      a}}{1-q}\right)^{-1/2}
= \hbar^{k/2} a^k q^{\frac{k^2}{2}} = 
\left(- \hbar^{1/2}\right)^k \frac{\vartheta(a)}{\vartheta(q^k a)}
\,.
$$
\end{proof}

\subsubsection{}
Let $\Stab(F)$ denote the restriction of elliptic stable envelopes
to a component $F\subset X^\bA$. 
By definition of stable envelopes, 
\begin{equation}
\left. \left(\frac{\Stab_\#(F)}
{\Theta(T^{1/2})}\right)_\tw \right/ \frac{\Stab_\#(F)}
{\Theta(T^{1/2})} 
= \dots \,\,  {z_\#}^{\bla(\cdot , \sigma)}  \,\label{stab_tw}
\end{equation}
where dots stand for a scalar factor that depends on the component 
$F$.  Degrees of the ``constant'' twisted quasimaps to $X^\sigma$ are 
computed as follows  
$$
\deg_\tw - \deg = - \bla(\cdot,\sigma) \,. 
$$
This gives the following 
\begin{Lemma}
The localization of 
$$
z^{\deg_\tw} \, \cK^{1/2}_{\vir,\tw} \,
\ev^*_\tw\left(\Stab_\#(F)\right) 
$$
depends on the twist only through a scalar factor that depends on
$F$. 
\end{Lemma}

\noindent 
Putting it all together, we obtain the following 

\begin{Proposition}\label{p_comp_tw} 
Under the assumption \eqref{tsigma}, 
the $\Ct_q$-localization of 
$$
z^{\deg_\tw} \, \cG \otimes \ev^*_\tw(\cF)\,, 
$$
for $\cG$ and $\cF$ as in \eqref{G1} and \eqref{F1}, depends on the 
twist only through a scalar function of the component of the fixed 
locus in the domain of stable envelope.  Therefore, the 
sum of such localization contributions for untwisted quasimaps 
is regular at \eqref{div4}. 
\end{Proposition}

\subsubsection{}

We now consider the general case, when the required twist of 
quasimaps to $X'$ cannot be obtained by twisting the quasimaps 
to the ambient $X$. 

By our assumption, $\bA$ is not in the kernel of $t$, therefore 
we can find $\Ct_a \subset \bA$ such that 
$$
t(a) = a^m\,, \quad m\ne 0 \,. 
$$
Restricted to $X'$, the tautological bundles of $\cV_i$ of $X$ 
will split according to characters of $\bT'$ 
$$
\cV_i \Big|_{X'} = \bigoplus_{\eta\in \left(\bT'\right)^\wedge}
\cV_{i,\eta} \,,
$$
and we can take a coarser decomposition
$$
\cV_i \Big|_{X'} = \bigoplus_{k=0}^{m-1} 
\cV_{i,k} 
$$
by the characters of 
$$
\Gamma= \Ct_a \cap \bT' \cong \Z/m \,. 
$$
The torus $\Ct_q$ has a well-defined action on 
\begin{equation}
\cV_{i,k,\textup{shift}} = q^{-k/m} \, \cV_{i,k}\label{cVshift}
\end{equation}
via the multivalued map 
$$
\sigma: \Ct_q \owns q \mapsto q^{1/m} \in \Ct_a \,. 
$$
Using this map, we can define the shifted and twisted 
tautological bundles over $\QM(X')_\tw$, and similarly for 
the framing bundles $\cW_{i,k,\textup{shift}}$. 

Note that shifts $q^{-k/m}$ for all these bundles lie in 
$(-1,0]$. Therefore, the shifts for the tangent bundle, which 
is a sesquilinear expression in $\cV_i$ and $\cW_i$, lie 
in $(-1,1)$ and the shift is zero precisely on 
$\Gamma$-invariants. This means that: 
\begin{itemize}
\item[---] shifts occur only in the virtual normal directions to
  $\QM(X')$ \,, 
\item[---] a polarization of
  $$
N^{1/2}=\left(T^{1/2}\big|_{X'}\right)_{\textup{$\bT'$ moving}}
$$ separates the shifts
into opposite pairs. 
\end{itemize}

We may consider $\QM(X')$ with a new obstruction 
theory, which is its own obstruction theory together with the 
contribution of shifted virtual normal bundle. 
This will define the sheaf $\tO_{\vir,\textup{shift}}$. 

\subsubsection{}

Consider the diagram of maps 
\begin{equation}
  \label{EllXX'}
\xymatrix{
  \Ell_\bT(X) \ar[d] & \Ell_\bT(X') \ar[l]_{\iota^*} \ar[r]^{t^*}
  \ar[d] & \Ell_{\Ct_t}(X') \ar[d] \\
  \cE_\bT & \cE_\bT \ar@{=}[l] \ar[r]^t& E 
}
\end{equation}
in which $\iota^*$ is the functorial map induced by the 
inclusion 
$$
\iota: X' \to X \,. 
$$
The shifts in \eqref{cVshift} induce an automorphism of 
$\Ell_\bT(X')$ which covers translation by $\sigma(q)$ on the 
base 
\begin{equation}
  \label{shiftX'}
\xymatrix{
  \Ell_\bT(X') \ar[d] \ar[rr]^{\textup{shift} } && \Ell_\bT(X') 
  \ar[d] \\
  \cE_\bT \ar[rr]^{t\mapsto \sigma(q) t} && \cE_\bT 
} \,. 
\end{equation}
We define 
\begin{align}
\cG & = \tO_{\vir,\textup{shift}} \otimes \textup{Tautological at $0$} \,, \label{G2} \\
\cF & = \left( \frac{\Stab_\# \otimes \left(\det T^{1/2}\right)^{-1/2}} 
{\Phi(q T^{\vee})} \right)_\textup{shift} \label{F2} \,,
\end{align}
where the shift of the stable envelope is the pullback under 
$(\iota \circ \textup{shift})$. 

After the shifts, the relation \eqref{stab_tw} is unchanged, as it 
concerns the degrees of curves in $X'$, which are not affected by a  shift in equivariant 
structure. The shift of polarization affects $\Theta(T^{1/2})$ in the 
obvious way. Namely, it becomes
$\Theta\left(T^{1/2}_\textup{shift}\right)$, 
and its transformation under $\sigma$ precisely matches the 
transformation of $\tO_{\vir,\textup{shift}}$. This is the computation 
we did in Sections \ref{s_twist_comp1} and \ref{s_twist_comp2}. 

This proves the conclusion of Proposition \ref{p_comp_tw}  without 
the assumption \eqref{tsigma} and concludes the proof of Theorem
\ref{t_ps}.

\end{document}